\documentclass[review, 3p, onecolumn]{elsarticle}

\usepackage[utf8x]{inputenc}
\usepackage{amsmath,amssymb,amsfonts}
\usepackage{graphicx}
\usepackage{bm}
\usepackage{booktabs}
\usepackage{siunitx}
\usepackage{xcolor}
\usepackage[hidelinks]{hyperref}
\usepackage{algorithm}
\usepackage{algpseudocode}
\usepackage{multirow}
\usepackage{lineno}
\usepackage{tikz}
\usetikzlibrary{fit,backgrounds}
\usetikzlibrary{arrows.meta,positioning,calc,decorations.pathmorphing,shadows}
\usetikzlibrary{decorations.pathreplacing}
\definecolor{cInput}{RGB}{30,100,220}     
\definecolor{cModel}{RGB}{120,60,190}     
\definecolor{cGate}{RGB}{235,90,20}       
\definecolor{cOut}{RGB}{20,130,60}        
\definecolor{cData}{RGB}{110,110,110}     

\newcommand{\Tmax}{T_{\max}}
\newcommand{\Nh}{N_h}

\journal{Computer Methods in Applied Mechanics and Engineering}

\begin{document}

\begin{frontmatter}

\title{
Model predictive control for laser thermal processing: operator learning, closed-loop validation, and out-of-distribution analysis}

\author[inst1]{Smajil Halilovic\corref{cor1}}
\ead{smajil.halilovic@tum.de}
\cortext[cor1]{Corresponding author}
\author[inst2]{Yuchen Zhang}
\author[inst2]{Daniel M. Tartakovsky}

\affiliation[inst1]{organization={Chair of Renewable and Sustainable Energy Systems, Technical University of Munich},
            city={Munich},
            country={Germany}} 
\affiliation[inst2]{organization={Department of Energy Science and Engineering, Stanford University},
            city={Stanford, CA}, country={USA}}

\begin{abstract}
Laser-based thermal processing, such as laser powder bed fusion, requires tight regulation of the peak surface temperature: heat accumulates where the moving source re-enters previously heated material, driving the temperature out of its process window and causing defects. High-fidelity  thermal models capture this physics but are too slow for online optimization, which motivates fast, differentiable, and generalizable surrogates. We develop and validate a complete surrogate-based control pipeline that regulates the maximum surface temperature of a moving laser  on a 304-stainless-steel substrate. We also determine conditions under which our surrogate can be trusted inside the control loop by probing its out-of-distribution limits. A key component of our surrogate
is a multi-step deep operator network bespoke for moving sources: its branch subnetwork encodes the future power and trajectory (position and velocity) sequence, while its trunk encodes the current peak temperature and the temperature at the future laser locations, yielding a one-shot five-step
prediction. By way of illustration, we use this surrogate as a smooth (algebraic-rectifier) nonlinear program inside a receding-horizon model predictive controller solved in CasADi/IPOPT. 
The surrogate forward pass is over thousand times faster than the equivalent finite-difference steps. We show that aggregate open-loop accuracy is necessary but not sufficient for control-readiness: two surrogates with near-identical offline error behave drastically differently in closed loop. A controlled two-ensemble data design reduces a 91~K path-corner underprediction
failure to 1.4~K, and a calibrated one-sided constraint margin of 13~K yields zero violations of the true upper bound on all tested paths.
\end{abstract}


\end{frontmatter}

\section{Introduction}
\label{sec:intro}
Laser-based manufacturing processes require precise regulation of the thermal field.
In laser powder bed fusion and directed energy deposition, part quality is governed by
keeping the melt-pool peak temperature within a narrow process window, and defects arise
when this window is violated~\cite{king2015lpbf,debroy2018am}. A dominant mechanism is
heat accumulation: when the moving source re-enters material it has recently heated, most
severely at sharp corners and path reversals, the peak temperature overshoots the upper
operating bound. Regulating the peak temperature along complex scan paths, using laser
power as the manipulated input, is therefore a central control problem, and the one addressed in this work.

High-fidelity transient conduction models resolve this physics but are far too slow to be
evaluated inside an optimizer, which must score many candidate power sequences per control
interval. This motivates surrogate models that are simultaneously fast, differentiable, and
able to generalize beyond their training distribution. Several families address this need.
Koopman and other lifting-based linear surrogates render the prediction linear in a lifted
state, so that model predictive control (MPC) reduces to a convex quadratic program and
generalization holds by construction~\cite{korda2018linear,zhang2022robusttube}; a regional
Koopman surrogate has recently been applied to the same distributed-thermal problem class as
studied here~\cite{zhang2026regional}. A second family embeds a neural network directly as the
prediction model inside the online optimization~\citep{salzmann2023realtimenmpc,karg2020efficient,hertneck2018learning,hewing2020learningbased}.
A third family, operator learning, approximates the solution operator of a parametric partial
differential equation; besides the Fourier neural operator~\cite{li2021fno}, the
deep operator network (DeepONet) rests on the operator universal-approximation theorem~\cite{chen1995universal,lu2021deeponet} and has been extended to full-field predictions in materials processing and additive manufacturing~\cite{he2024sdeeponet,kushwaha2024advanced,olearyroseberry2024dino}.
The multi-step DeepONet (MS-DeepONet) of de Jong et al.~\cite{dejong2025deeponetmpc} is the
enabling architecture that makes non-recurrent one-shot horizon prediction compatible with an
optimization graph. In additive manufacturing specifically, learned thermal control has been
pursued with deep reinforcement learning~\cite{ogoke2021thermal} and, most closely to the
present work, a one-shot multi-step time-series surrogate driving a real-time MPC of directed
energy deposition~\cite{chen2025tidempc}.

Despite this activity, no prior work embeds an operator network inside an online constrained nonlinear-program MPC for a moving-source distributed thermal system with a peak-temperature constraint, together with an analysis of its out-of-distribution (OOD) behavior.
The objective of this paper is therefore twofold. 
First, we develop and validate a complete surrogate-based control pipeline: from excitation-aware training-data generation, through a moving-source-aware operator-network architecture, to a differentiable embedding in a receding-horizon controller, that regulates the peak surface temperature of a moving laser source within a prescribed process window along complex, fixed scan paths, using only measurements available from a thermal camera. 
Second, we ask under what conditions such a learned surrogate can be trusted inside the control loop: we deliberately expose its out-of-distribution limits on geometrically unseen scan paths, diagnose the resulting failure mode, and demonstrate how it is mitigated by targeted training-data design and a calibrated constraint margin.
In pursuit of this objective, we make five contributions. (C1) A moving-source-aware MS-DeepONet branch/trunk design in which the trunk ingests temperatures at the future laser locations, encoding pre-heated zones. (C2)~A causality-preserving, excitation-aware
two-ensemble methodology for generating training data. (C3)~A differentiable smooth-rectifier
embedding of the surrogate as a nonlinear program, with a receding-horizon soft-constrained
controller. (C4)~A systematic OOD diagnosis, a $91$~K corner under-prediction at a sharp laser path
reversal, together with a data-design remedy and a calibrated one-sided constraint
tightening ($\delta=13$~K) that yields zero violations of the true upper bound. (C5)~Closed-loop
validation against a high-fidelity finite-difference plant emulating thermal-camera feedback,
with an objective computational and real-time analysis. 
A cross-cutting methodological finding, answering the second objective, is that aggregate open-loop validation is insufficient to certify a learned surrogate for control,
because the failure regime is rarely sampled by open-loop excitation yet is permanently
occupied by the controller. The remainder of the paper is organized as follows.
Section~\ref{sec:problem} states the thermal control problem and the high-fidelity solver;
Section~\ref{sec:surrogate} develops the moving-source-aware MS-DeepONet and the data design;
Section~\ref{sec:mpc} presents the differentiable NLP embedding and the receding-horizon
controller; Section~\ref{sec:results} reports offline and closed-loop results, the OOD
diagnosis, and ablations; Section~\ref{sec:discussion} discusses the findings; and
Section~\ref{sec:conclusion} concludes.

\section{Problem formulation}
\label{sec:problem}

\subsection{Governing equations and physical setup}
\label{sec:physics}

\begin{figure}[t]
    \centering
    \begin{tikzpicture}[
        >={Latex[length=2mm]},
        font=\small,
        box/.style={rounded corners=2pt}
    ]
        \coordinate (A) at (0,0);
        \draw[fill=cInput!8, draw=cInput!60!black, thick]
            (0,0) -- (6.0,0) -- (7.4,0.9) -- (1.4,0.9) -- cycle;      
        \draw[fill=cInput!15, draw=cInput!60!black, thick]
            (0,0) -- (6.0,0) -- (6.0,-1.1) -- (0,-1.1) -- cycle;      
        \draw[fill=cInput!20, draw=cInput!60!black, thick]
            (6.0,0) -- (7.4,0.9) -- (7.4,-0.2) -- (6.0,-1.1) -- cycle;

        \draw[cModel, very thick, ->, rounded corners=1pt]
            (1.6,0.30) -- (5.2,0.30) -- (2.6,0.62) -- (6.2,0.62);
        \node[cModel, anchor=west, font=\scriptsize] at (5.7,0.30)
            {$\bm{\xi}(t)$};

        \coordinate (spot) at (4.1,0.46);
        \fill[cGate!80, opacity=0.25] (spot) ellipse (0.32 and 0.10);
        \fill[cGate]                  (spot) ellipse (0.10 and 0.035);
        \draw[cGate, thick, fill=cGate!20, opacity=0.65]
            ($(spot)+(0,2.0)$) -- ($(spot)+(-0.30,0.05)$)
            arc[start angle=180, end angle=360, x radius=0.30, y radius=0.09]
            -- cycle;
        \node[cGate!80!black, anchor=west, font=\scriptsize]
            at ($(spot)+(0.42,1.55)$) {$P(t),\ \sigma$};

        \draw[cOut, thick, ->] ($(spot)+(-1.9,1.15)$) -- ($(spot)+(-0.35,0.12)$);
        \node[cOut, anchor=east, font=\scriptsize] at ($(spot)+(-1.9,1.25)$)
            {$\Tmax(t)=\max_{(x,y)} T(x,y,L_z,t)$};

        \node[cData, font=\scriptsize, anchor=west] at (0.05,-1.35)
            {$T = T_0$ (bottom)};
        \node[cData, font=\scriptsize, anchor=west] at (3.6,-1.35)
            {adiabatic lateral walls};
        \node[font=\scriptsize, anchor=west] at (0.1,1.15)
            {$\Omega:\ 15 \times 10 \times 2\,$mm};
    \end{tikzpicture}
    \caption{Physical setup: a Gaussian laser source of power $P(t)$ and
    width $\sigma$ traverses the prescribed scan path $\bm{\xi}(t)$ on the
    top surface of the substrate. The controlled output is the maximum
    surface temperature $\Tmax(t)$.}
    \label{fig:setup}
\end{figure}
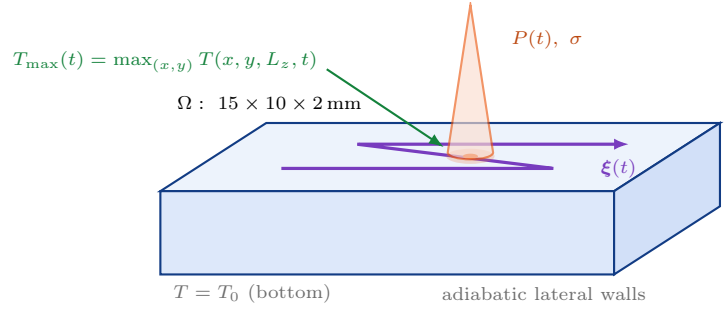

\begin{figure}[b]
    \centering
    \includegraphics[width=0.95\textwidth]{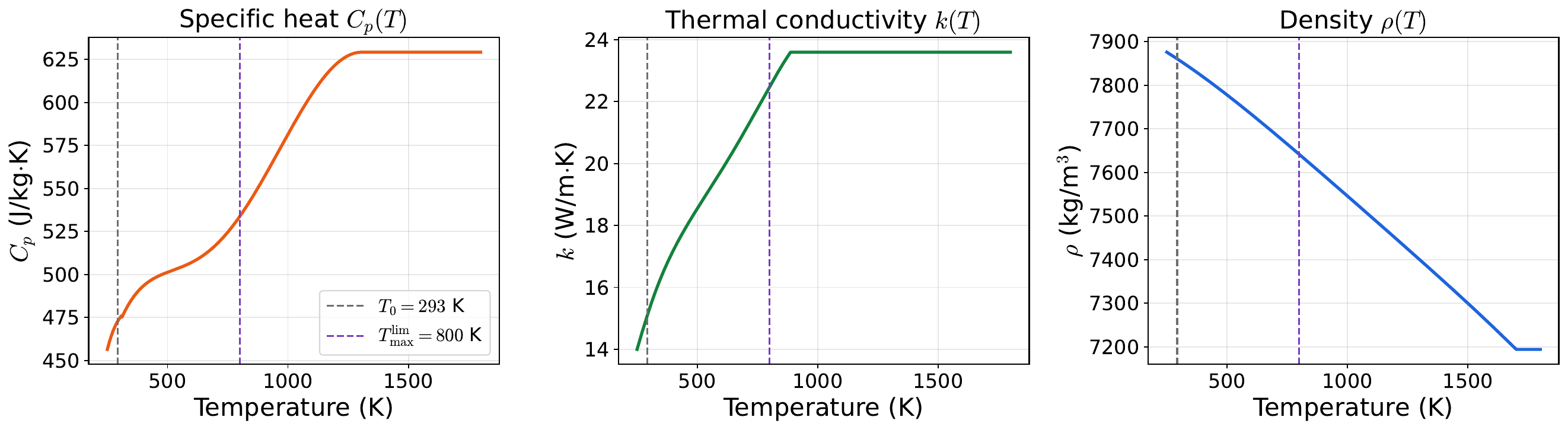}
    \caption{Temperature-dependent density $\rho(T)$, specific heat $c_p(T)$, and thermal conductivity $k(T)$ of 304 stainless steel, from~\cite{kim1975stainless,mills2002thermophysical}.}
    \label{fig:material}
\end{figure}

We adopt the laser thermal-control problem of Zhang and
Tartakovsky~\cite{zhang2026regional}, illustrated in
Fig.~\ref{fig:setup}. 
We consider transient heat conduction in a metallic workpiece occupying the
cuboid $\Omega = [-L_x/2, L_x/2] \times [-L_y/2, L_y/2] \times [0, L_z]
\subset \mathbb{R}^3$, with in-plane coordinates $\mathbf{x}_s = (x, y)^\top$
and depth coordinate $z$ measured from the bottom surface. The temperature
field $T(\mathbf{x}, t)$, $\mathbf{x} = (x, y, z)^\top$, evolves according to
the nonlinear heat equation
\begin{equation}
    \rho(T)\, c_p(T)\, \frac{\partial T}{\partial t}
    \;=\;
    \nabla \cdot \bigl( k(T)\, \nabla T \bigr)
    \qquad \text{in } \Omega \times \mathbb{R}^+,
    \label{eq:heat}
\end{equation}
where $\rho$ [kg\,m$^{-3}$], $c_p$ [J\,kg$^{-1}$K$^{-1}$], and $k$
[W\,m$^{-1}$K$^{-1}$] denote the temperature-dependent density, specific
heat, and thermal conductivity of 304 stainless steel; their dependence on
$T$ is shown in Fig.~\ref{fig:material}. The moving laser is modeled as a
Gaussian heat flux applied to the top surface $z = L_z$,
\begin{equation}
    -k(T)\, \nabla T \cdot \mathbf{n}
    \;=\;
    \frac{P(t)}{2\pi\sigma^2}
    \exp\!\left(
        -\frac{\bigl\| \mathbf{x}_s - \bm{\xi}(t) \bigr\|^2}{2\sigma^2}
    \right),
    \qquad z = L_z,
    \label{eq:laser}
\end{equation}
where $P(t) \in [0, P_{\max}]$ [W] is the laser power (the control input),
$\sigma = 100\,\mu$m the beam width, and
$\bm{\xi}(t) = \bigl( X_c(t), Y_c(t) \bigr)^\top$ the prescribed scan path.
The lateral boundaries are adiabatic, the bottom surface is held at the
ambient temperature $T_0$, and the initial condition is uniform,
$T(\mathbf{x}, 0) = T_0$:
\begin{equation}
    -k(T)\, \nabla T \cdot \mathbf{n} = 0 \;\; \text{on } \partial\Omega_{\text{side}},
    \qquad
    T = T_0 \;\; \text{on } \partial\Omega_{\text{bottom}}.
    \label{eq:bcs}
\end{equation}

The quantity to be regulated is the maximum surface temperature,
\begin{equation}
    \Tmax(t)
    \;=\;
    \max_{(x,y)} \; T\bigl( x, y, L_z, t \bigr),
    \label{eq:tmax}
\end{equation}
which governs melt-pool stability and defect formation. The control
objective is to keep $\Tmax$ within the prescribed process window
\begin{equation}
    \Tmax(t) \;\in\; [\,T_{\min},\, T_{\max}^{\lim}\,]
    \;=\; [\,760,\, 800\,]~\text{K}
    \qquad \forall t,
    \label{eq:window}
\end{equation}
by modulating $P(t)$ along a scan path $\bm{\xi}(t)$ that is fixed by the
part geometry and therefore not a degree of freedom of the controller.

\subsection{High-fidelity finite-difference model (digital ground truth)}
\label{sec:fdm}

Equations~\eqref{eq:heat}--\eqref{eq:bcs} are discretized on a structured
Cartesian grid of $N_x \times N_y \times N_z = 151 \times 101 \times 21$
nodes with central differences in space and an implicit time-marching
scheme. 
The nonlinear diffusion term is resolved by Newton iterations at
each step. A fixed time step $\Delta t = 1.25 \times 10^{-4}$~s is used, and
a single production run spans $N_t = 401$ states ($t_f = 0.05$~s). We refer
to this solver as the finite-difference model (FDM) and treat its
(noise-free) solution as the physics-based ground truth.

The FDM serves two distinct roles. First, it is the \emph{training-data
generator}: ensembles of simulations driven by randomized power profiles and
scan-path geometries (Section~\ref{sec:data}) provide the input--output
pairs from which the surrogate is learned. Second, it is the
\emph{in-the-loop plant} during closed-loop validation: at every control
step the FDM advances the true thermal state, and its top-surface field
emulates the thermal-camera feedback (the measured $\Tmax$ and the local
surface temperatures at future laser positions) supplied to the controller.

The FDM itself is unsuitable for control. 
A single high-fidelity run (401 steps) requires on the order of a minute of
wall-clock time on the benchmark machine ($\sim\!103$~ms per implicit step,
$\sim\!41$~s per run), whereas the controller must evaluate candidate power
sequences many times per control interval, which rules out direct in-the-loop use.
This gap of several orders of magnitude motivates replacing the FDM inside
the optimization loop by a fast learned surrogate, while retaining the FDM
as the external referee that validates the resulting closed-loop behavior.

\subsection{Model predictive control problem}
\label{sec:mpc}

Following the receding-horizon paradigm~\cite{mayne2014mpc,rawlings2017mpc},
at discrete time $t_k = k \Delta t$ the controller selects the power
sequence $\mathbf{P}_k = ( P_{k+1}, \dots, P_{k+\Nh} )^\top$ over a horizon
of $\Nh$ steps by solving the nonlinear program
\begin{subequations}
\label{eq:mpc}
\begin{align}
    \min_{\mathbf{P}_k,\, \bm{\epsilon}_k} \quad
    & w_P \sum_{i=1}^{\Nh} \bar{P}_{k+i}^{\,2}
    \;+\; w_{\Delta P} \sum_{i=1}^{\Nh}
        \bigl( \bar{P}_{k+i} - \bar{P}_{k+i-1} \bigr)^2
    \;+\; w_\epsilon \sum_{i=1}^{\Nh} \epsilon_{k+i}
    \label{eq:mpc-cost} \\[2pt]
    \text{s.t.} \quad
    & \widehat{T}_{\max,k+i}
        \;=\; \bigl[ \mathcal{F}_{\text{surr}}
        ( \mathbf{u}_k, \mathbf{y}_k ) \bigr]_i,
    && i = 1, \dots, \Nh,
    \label{eq:mpc-dyn} \\
    & \widehat{T}_{\max,k+i} \;\le\; T_{\max}^{\lim},
    && i = 1, \dots, \Nh,
    \label{eq:mpc-hard} \\
    & \widehat{T}_{\max,k+i} + \epsilon_{k+i} \;\ge\; T_{\min},
    \quad \epsilon_{k+i} \ge 0,
    && i = 1, \dots, \Nh,
    \label{eq:mpc-soft} \\
    & 0 \;\le\; P_{k+i} \;\le\; P_{\max},
    && i = 1, \dots, \Nh,
    \label{eq:mpc-box}
\end{align}
\end{subequations}
where $\bar{P} = P / P_{\max}$ denotes normalized power,
$\bar{P}_k$ is the previously applied input (so that the rate penalty spans
the receding-horizon boundary), and
$\bm{\epsilon}_k = ( \epsilon_{k+1}, \dots, \epsilon_{k+\Nh} )^\top$ are
slack variables. The upper bound~\eqref{eq:mpc-hard} is enforced as a hard
constraint, since exceeding $T_{\max}^{\lim}$ risks material damage, whereas
the lower bound~\eqref{eq:mpc-soft} is softened by an exact
$\ell_1$ penalty with $w_\epsilon \gg w_P, w_{\Delta P}$: whenever the lower
bound is physically attainable, the optimizer drives
$\bm{\epsilon}_k \to \mathbf{0}$ and the constraint is satisfied exactly,
while infeasible transients (e.g., during initial heat-up) do not render the
program infeasible. We use $w_P = 1$, $w_{\Delta P} = 10$,
$w_\epsilon = 10^6$, and $P_{\max} = 20$~W.

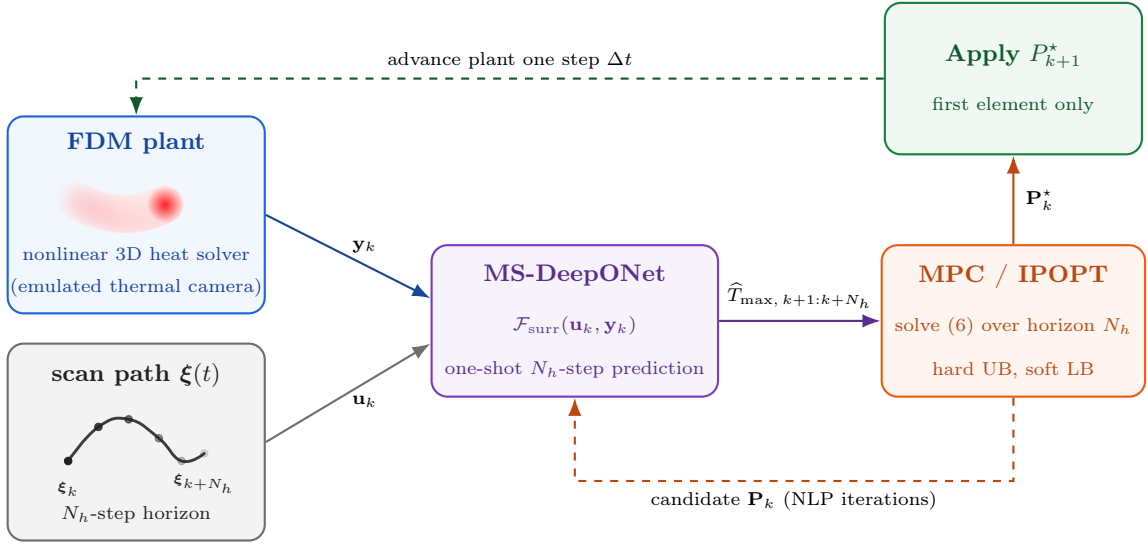
\begin{figure}[t]
    \centering
    \begin{tikzpicture}[
        >={Latex[length=2.4mm]},
        font=\small,
        blk/.style={rounded corners=6pt, thick, align=center,
                    inner sep=5pt, minimum width=34mm, minimum height=26mm},
        plantb/.style={blk, draw=cInput, fill=cInput!6, text=cInput!70!black},
        surrb/.style ={blk, draw=cModel, fill=cModel!6, text=cModel!75!black,
                       minimum height=20mm},
        optb/.style  ={blk, draw=cGate,  fill=cGate!6,  text=cGate!80!black,
                       minimum height=20mm},
        outb/.style  ={blk, draw=cOut,   fill=cOut!6,   text=cOut!70!black,
                       minimum height=20mm},
        datab/.style ={blk, draw=cData,  fill=cData!8,  text=cData!30!black},
        cameo/.style ={inner sep=0pt, outer sep=0pt},
        lbl/.style   ={font=\scriptsize, align=center, text=black},
        flow/.style  ={->, thick},
    ]
    \def\dx{58mm}                       
    \coordinate (Pc) at (0,0);          
    \coordinate (Sc) at (\dx,0);        
    \coordinate (Mc) at (2*\dx,0);      
    \coordinate (Ac) at (2*\dx,32mm);   
 
    \node[plantb] (plant) at ($(Pc)+(0,14mm)$) {};   
    \node[lbl, anchor=north, text=cInput!70!black, font=\small\bfseries]
        at ($(plant.north)+(0,-1.0mm)$) {FDM plant};
    \begin{scope}
      \clip ($(plant.center)+(-15mm,-8.5mm)$) rectangle ($(plant.center)+(15mm,9mm)$);
    
      \coordinate (Ph) at ($(plant.center)+(4mm,1.5mm)$);   
      \coordinate (Bc) at ($(Ph)+(-6mm,-3mm)$);             
      \coordinate (Pt) at ($(Ph)+(-12mm,0.5mm)$);           
    
      \def\rad{2.5mm}                                        
    
      \foreach \t in {0,0.02,...,1}{
        \pgfmathsetmacro{\op}{0.12 + 0.83*\t*\t}            
        \pgfmathtruncatemacro{\lvl}{30 + 60*\t}             
        \coordinate (Q) at
          ($($(Pt)!\t!(Bc)$)!\t!($(Bc)!\t!(Ph)$)$);         
        \shade[inner color=red!\lvl, outer color=red!8, opacity=\op]
          (Q) circle (\rad);
      }
    \end{scope}
    \node[lbl, anchor=south, text=cInput!70!black]
        at ($(plant.south)+(0,1.15mm)$)
        {\scriptsize nonlinear 3D heat solver\\\scriptsize (emulated thermal camera)};
 
    \node[surrb] (surr) at (Sc)
        {\textbf{MS-DeepONet}\\[2pt]
         \scriptsize $\mathcal{F}_{\text{surr}}(\mathbf{u}_k,\mathbf{y}_k)$\\
         \scriptsize one-shot $\Nh$-step prediction};
    \node[optb] (mpc) at (Mc)
        {\textbf{MPC / IPOPT}\\[2pt]
         \scriptsize solve \eqref{eq:mpc} over horizon $\Nh$\\
         \scriptsize hard UB, soft LB};
    \node[outb] (act) at (Ac)
        {\textbf{Apply} $P_{k+1}^{\star}$\\[2pt]
         \scriptsize first element only};
 
    \node[datab] (path) at ($(Pc)+(0,-16mm)$) {};     
    \node[lbl, anchor=north, text=cData!30!black, font=\small\bfseries]
        at ($(path.north)+(0,-1.5mm)$) {scan path $\bm{\xi}(t)$};
    \node[lbl, anchor=south, text=cData!30!black]
        at ($(path.south)+(0,1.2mm)$) {\scriptsize $\Nh$-step horizon};

    \begin{scope}[shift={($(path.center)+(0,0.5mm)$)}]
        \draw[cData!55!black, line width=1.1pt, line cap=round]
            plot[smooth, tension=0.8] coordinates
            {(-9mm,-3mm) (-5mm,1.5mm) (-1mm,2.5mm) (3mm,0mm) (6mm,-3mm) (9mm,-2mm)};
        \foreach \p/\op in {%
            {-9mm,-3mm}/1.0, {-5mm,1.5mm}/0.82, {-1mm,2.5mm}/0.64,
            {3mm,0mm}/0.46, {6mm,-3mm}/0.28, {9mm,-2mm}/0.14}{
            \fill[cData!30!black, opacity=\op] (\p) circle (1.4pt);
            \draw[cData!30!black, opacity=\op, line width=0.4pt] (\p) circle (1.4pt);
        }
        \node[font=\tiny, text=cData!30!black, anchor=north]
            at (-9mm,-4.5mm) {$\bm{\xi}_k$};
        \node[font=\tiny, text=cData!30!black, anchor=north]
            at (9mm,-3.5mm) {$\bm{\xi}_{k+\Nh}$};
    \end{scope}
 
    \draw[flow, cInput!70!black]
        (plant.east) -- node[lbl, above=1pt, pos=0.6] {$\mathbf{y}_k$}
        ($(surr.west)+(0,3mm)$);
    \draw[flow, cModel!75!black]
        (surr) -- node[lbl, above=1pt, pos=0.5]
        {$\widehat{T}_{\max,\,k+1:k+\Nh}$} (mpc);
    \draw[flow, cGate!80!black]
        (mpc) -- node[lbl, right=1pt] {$\mathbf{P}_k^{\star}$} (act);
 
    \draw[flow, cGate!80!black, dashed]
        (mpc.south) -- ++(0,-11mm)
        -| node[lbl, below, pos=0.25]
           {candidate $\mathbf{P}_k$ (NLP iterations)}
        (surr.south);
 
    \draw[flow, cOut!70!black, dashed]
        (act.west) -| (plant.north);
    \node[lbl, above, text=black]
        at ($(plant.north)!0.5!(plant.north-|act.west) + (0,5mm)$)
        {advance plant one step $\Delta t$};
 
    \draw[flow, cData]
        (path.east) -- node[lbl, below=1pt, pos=0.6] {$\mathbf{u}_k$}
        ($(surr.west)+(0,-3mm)$);
 
    \end{tikzpicture}
    \caption{Closed-loop architecture. At each control step the FDM plant,
    which emulates a thermal camera, supplies the measured thermal state
    $\mathbf{y}_k$ (the peak temperature and the surface temperatures at the
    upcoming laser locations). The MS-DeepONet predicts the full
    $\Nh$-step $\Tmax$ trajectory in a single evaluation, queried repeatedly
    by IPOPT as it solves~\eqref{eq:mpc}; the first optimal power
    $P_{k+1}^{\star}$ is applied, the plant advances by one step $\Delta t$,
    and the horizon recedes.}
    \label{fig:closedloop}
\end{figure}

The dynamics constraint~\eqref{eq:mpc-dyn} is supplied by the learned
surrogate $\mathcal{F}_{\text{surr}}$, which maps the planned actuation
and kinematics over the horizon, $\mathbf{u}_k$, together with the measured thermal state $\mathbf{y}_k$ to the predicted trajectory
$( \widehat{T}_{\max,k+1}, \dots, \widehat{T}_{\max,k+\Nh} )^\top$ in a
single evaluation; its construction is detailed in
Section~\ref{sec:surrogate}. Following the receding-horizon principle, only
the first element $P_{k+1}^\star$ of the optimal sequence is applied to the
plant; the state is then re-measured and problem~\eqref{eq:mpc} is solved
anew at $t_{k+1}$. The closed loop is summarized in
Fig.~\ref{fig:closedloop}.

\section{Multi-step DeepONet surrogate for moving-source thermal dynamics}
\label{sec:surrogate}

\subsection{MS-DeepONet architecture}
\label{sec:arch}

The surrogate $\mathcal{F}_{\text{surr}}$ in~\eqref{eq:mpc-dyn} is a
multi-step deep operator network (MS-DeepONet) that predicts the entire
$\Nh$-step $\Tmax$ trajectory in a single forward pass. It comprises two
feedforward subnetworks: a \emph{branch} network encoding the planned
actuation and scan kinematics over the horizon, and a \emph{trunk} network
encoding the current thermal state. The branch maps its input
$\mathbf{u}_k$ to a coefficient matrix
$\mathbf{B}(\mathbf{u}_k) \in \mathbb{R}^{\Nh \times p}$ and the trunk maps
$\mathbf{y}_k$ to a basis vector
$\mathbf{T}(\mathbf{y}_k) \in \mathbb{R}^{p}$; the prediction is their inner
product per horizon step,
\begin{equation}
    \widehat{T}_{\max,k+i}
    \;=\;
    \sum_{j=1}^{p} B_{ij}(\mathbf{u}_k)\, T_{j}(\mathbf{y}_k) \;+\; b_i,
    \qquad i = 1, \dots, \Nh,
    \label{eq:msdeeponet}
\end{equation}
with a trainable bias $\mathbf{b} \in \mathbb{R}^{\Nh}$. Both subnetworks
use two hidden layers of $128$ neurons with rectified-linear (ReLU)
activations~\cite{nair2010relu}, and we employ
$p = 100$ basis functions; the sensitivity of the accuracy to $p$ and to the
layer width is quantified in Section~\ref{sec:ablations}. This deliberately
compact architecture keeps the nonlinear program~\eqref{eq:mpc} tractable:
the surrogate is evaluated and differentiated at every solver iteration, so
its cost enters the per-step control budget multiplicatively.

Two properties make this operator formulation suited to receding-horizon
control. First, DeepONets inherit a universal-approximation guarantee for
nonlinear operators~\cite{chen1995universal,lu2021deeponet,dejong2025deeponetmpc}, so a single network can represent the map from actuation sequences and thermal states to
temperature trajectories. Second, and central to this work, the prediction
in~\eqref{eq:msdeeponet} is \emph{one-shot}: the horizon is produced by a
single evaluation, with no autoregressive rollout. Recursive multi-step
predictors compose their own outputs $\Nh$ times, which is known to compound
single-step errors along the horizon~\cite{bentaieb2016bias}; the
one-shot structure avoids this mechanism by construction.
Section~\ref{sec:ablations} quantifies the resulting accuracy gap against
an autoregressive alternative trained on identical data.

\subsection{Moving-source-aware input design}
\label{sec:inputs}

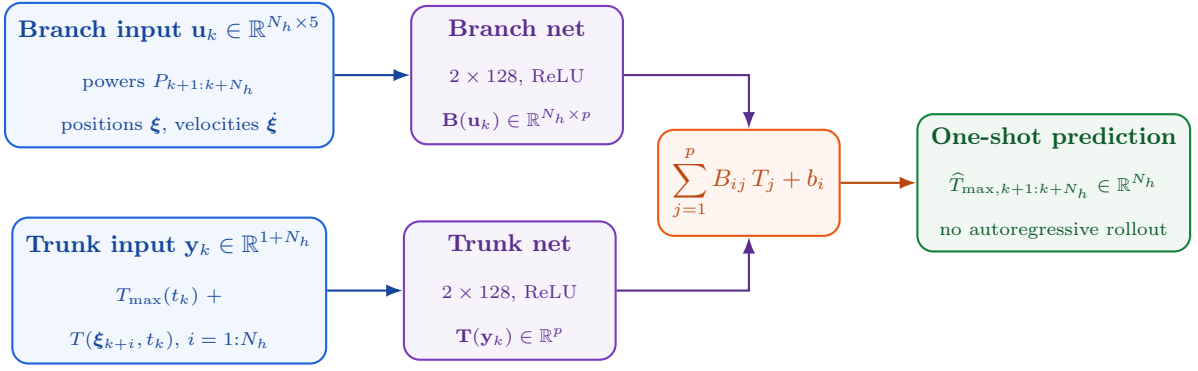
\begin{figure}[t]
    \centering
    \begin{tikzpicture}[
        >={Latex[length=2.2mm]},
        font=\small,
        blk/.style={rounded corners=6pt, thick, align=center, inner sep=5pt},
        inpb/.style ={blk, draw=cInput, fill=cInput!6, text=cInput!70!black,
                      minimum height=16mm, minimum width=36mm},
        netb/.style ={blk, draw=cModel, fill=cModel!6, text=cModel!75!black,
                      minimum height=16mm, minimum width=28mm},
        opb/.style  ={blk, draw=cGate,  fill=cGate!6,  text=cGate!80!black,
                      minimum height=14mm, minimum width=24mm},
        outbb/.style={blk, draw=cOut,   fill=cOut!6,   text=cOut!70!black,
                      minimum height=14mm, minimum width=32mm},
        lbl/.style  ={font=\scriptsize, align=center}
    ]
        \node[inpb] (u)
            {\textbf{Branch input} $\mathbf{u}_k \in \mathbb{R}^{\Nh\times 5}$\\[2pt]
             \scriptsize powers $P_{k+1:k+\Nh}$\\
             \scriptsize positions $\bm{\xi}$, velocities $\dot{\bm{\xi}}$};
        \node[inpb, below=9mm of u] (y)
            {\textbf{Trunk input} $\mathbf{y}_k \in \mathbb{R}^{1+\Nh}$\\[2pt]
             \scriptsize $\Tmax(t_k)$ $+$ \\
             \scriptsize $T(\bm{\xi}_{k+i}, t_k)$, $i=1{:}\Nh$};

        \node[netb, right=10mm of u] (bnet)
            {\textbf{Branch net}\\[1pt]
             \scriptsize $2\times128$, ReLU\\
             \scriptsize $\mathbf{B}(\mathbf{u}_k)\in\mathbb{R}^{\Nh\times p}$};
        \node[netb, right=10mm of y] (tnet)
            {\textbf{Trunk net}\\[1pt]
             \scriptsize $2\times128$, ReLU\\
             \scriptsize $\mathbf{T}(\mathbf{y}_k)\in\mathbb{R}^{p}$};

        \node[opb] (dot) at ($(bnet.south east)!0.5!(tnet.north east)+(17mm,0)$)
            {$\displaystyle \sum_{j=1}^{p} B_{ij}\, T_j + b_i$};

        \node[outbb, right=10mm of dot] (out)
            {\textbf{One-shot prediction}\\[1pt]
             \scriptsize $\widehat{T}_{\max,k+1:k+\Nh}\in\mathbb{R}^{\Nh}$\\
             \scriptsize no autoregressive rollout};

        \draw[->, thick, cInput!70!black] (u)  -- (bnet);
        \draw[->, thick, cInput!70!black] (y)  -- (tnet);
        \draw[->, thick, cModel!75!black] (bnet.east) -| ($(dot.north)+(0,0)$);
        \draw[->, thick, cModel!75!black] (tnet.east) -| ($(dot.south)+(0,0)$);
        \draw[->, thick, cGate!80!black]  (dot) -- (out);
    \end{tikzpicture}
    \caption{MS-DeepONet surrogate. The branch network encodes the causally
    aligned actuation and scan kinematics over the horizon
    (Eq.~\eqref{eq:branch}); the trunk network encodes the current thermal
    state including the spatial pre-heat preview at future laser locations
    (Eq.~\eqref{eq:trunk}); their inner product with bias yields the full
    $\Nh$-step prediction in a single evaluation.}
    \label{fig:architecture}
\end{figure}

The input construction, summarized in Fig.~\ref{fig:architecture}, encodes
two physical facts about moving-source heating: (i) the temperature response
over $[t_k, t_{k+\Nh}]$ is driven by the power \emph{applied during} that
interval along the \emph{path traversed} in it, and (ii) $\Tmax$ closely
tracks the instantaneous laser position except where the path revisits
recently heated material, in which case local pre-heating dominates.

\paragraph{Branch: causally aligned actuation and kinematics}
The power $P_{k+1}$ acts on the interval $[t_k, t_{k+1}]$ and therefore
drives the transition to the state at $t_{k+1}$. Respecting this causal
alignment, the $i$-th row of the branch input pairs the power of the
\emph{upcoming} interval with the position and velocity at its start:
\begin{equation}
    \mathbf{u}_k
    \;=\;
    \begin{bmatrix}
        P_{k+1}   & X_{c,k}      & Y_{c,k}      & \dot{X}_{c,k}      & \dot{Y}_{c,k} \\
        P_{k+2}   & X_{c,k+1}    & Y_{c,k+1}    & \dot{X}_{c,k+1}    & \dot{Y}_{c,k+1} \\
        \vdots    & \vdots       & \vdots       & \vdots            & \vdots \\
        P_{k+\Nh} & X_{c,k+\Nh-1} & Y_{c,k+\Nh-1} & \dot{X}_{c,k+\Nh-1} & \dot{Y}_{c,k+\Nh-1}
    \end{bmatrix}
    \in \mathbb{R}^{\Nh \times 5},
    \label{eq:branch}
\end{equation}
where $\bm{\xi}_j = \bm{\xi}(t_j) = (X_{c,j}, Y_{c,j})^\top$ is the laser position and
$\dot{\bm{\xi}}_j = \bigl( \bm{\xi}_{j+1} - \bm{\xi}_j \bigr)/\Delta t
= (\dot{X}_{c,j}, \dot{Y}_{c,j})^\top$ its finite-difference velocity.
The
velocities are not redundant: they encode the local dwell time of the
source, which sets the energy deposited per unit path length -- sharp turns
and slow traversals concentrate heat even at constant power.

\paragraph{Trunk: current state and spatial pre-heat preview}
The trunk observes the current maximum temperature together with the
temperatures of the material at the \emph{future} laser locations, sampled
from the current field:
\begin{equation}
    \mathbf{y}_k
    \;=\;
    \Bigl(
        \Tmax(t_k),\;
        T\bigl( X_{c,k+1}, Y_{c,k+1}, t_k \bigr),\;
        \dots,\;
        T\bigl( X_{c,k+\Nh}, Y_{c,k+\Nh}, t_k \bigr)
    \Bigr)^{\!\top}
    \in \mathbb{R}^{1+\Nh}.
    \label{eq:trunk}
\end{equation}
The look-ahead entries act as a spatial pre-heat preview: when the path
approaches previously heated material (most prominently at sharp corners,
where the source re-enters its own thermal wake) these temperatures rise
ahead of the laser and inform the network of the impending accumulation. In
closed loop they are read directly from the emulated thermal-camera image at
the (known) future path locations, so~\eqref{eq:trunk} requires no
quantities beyond those available to the controller.
An alternative trunk encoding that appends the coordinates of
the current $\Tmax$ location, $\bigl( X_{\max}(t_k), Y_{\max}(t_k) \bigr)$,
is evaluated in the ablation of Section~\ref{sec:ablations}.

\subsection{Excitation-aware training-data generation}
\label{sec:data}
 
Because the surrogate replaces the plant model inside the optimizer, its
training distribution must cover the actuation frequencies and path
kinematics that the closed loop will visit. To isolate the effect of the
path-geometry coverage on out-of-distribution behavior, we construct two
training ensembles that are identical in every respect except their geometry composition, and train one surrogate on each:
\begin{itemize}
    \item \textbf{Baseline ensemble} $\mathcal{D}_{\mathrm{base}}$:
    320 FDM simulations (of $N_t = 401$ states each) whose scan paths are drawn exclusively from the
    three "smooth" geometry classes: horizontal raster, vertical raster,
    and Archimedean spiral, in equal shares.
    \item \textbf{Corner-augmented ensemble} $\mathcal{D}_{\mathrm{corner}}$:
    320 FDM simulations in which the same three classes are complemented
    by a corner-rich random-polyline class introduced below, in the
    ratio 20:20:20:40.
\end{itemize}
The surrogate trained on $\mathcal{D}_{\mathrm{base}}$ serves as the
reference against which we expose the corner failure mode
(Section~\ref{sec:offline-results}); the surrogate trained on
$\mathcal{D}_{\mathrm{corner}}$ quantifies both the benefit of the
augmentation on unseen sharp-cornered geometry and its cost on
the seen geometries. We emphasize that this cost is measured, not
assumed away: rebalancing a fixed-size training budget toward
under-represented regimes necessarily dilutes the density of
in-distribution samples, and the resulting trade-off is reported in
Section~\ref{sec:offline-results}.

\begin{figure}[tb]
    \centering
    \includegraphics[width=0.9\linewidth]{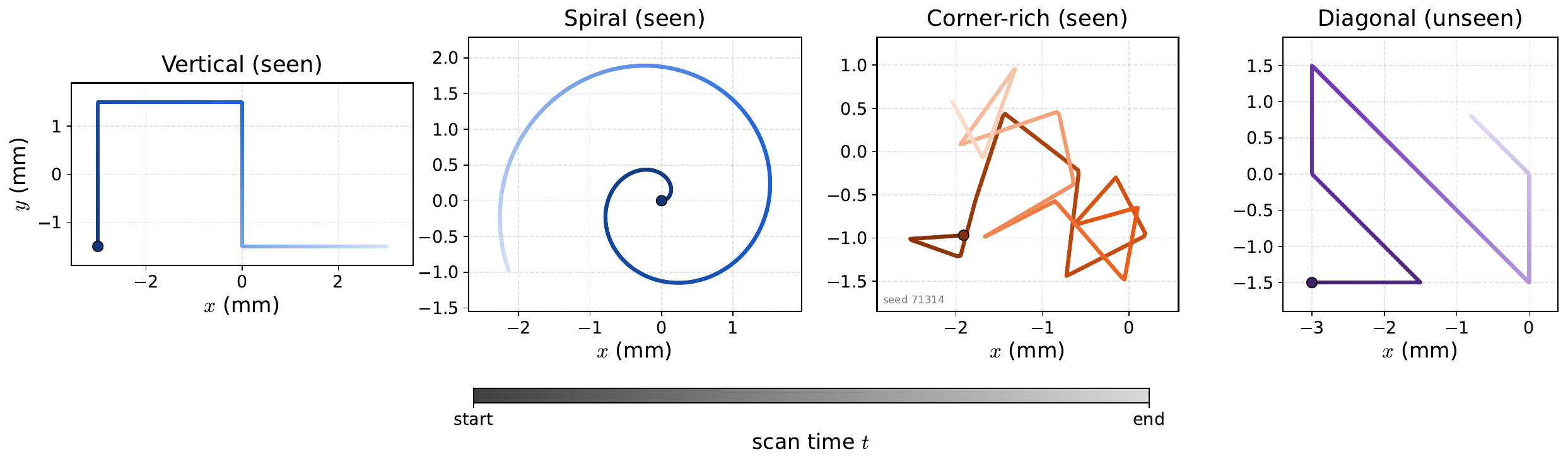}
    \caption{Path geometries: seen (vertical, spiral, example of random-polyline) and unseen (diagonal).}
    \label{fig:path-profiles}
\end{figure}

\paragraph{Path-geometry library}
The three "smooth" classes provide sustained straight traverses and mostly smooth
curvature; on their own, however, they under-represent sharp direction reversals, whose thermal signature (the source re-entering its own wake) is precisely the regime in which a surrogate trained on smooth
paths fails (Section~\ref{sec:offline-results}). Each corner-rich path
of $\mathcal{D}_{\mathrm{corner}}$ is a random polyline with segment
lengths drawn from $[0.12, 0.25]$~mm, corresponding to the distance
$v N_h \Delta t \approx 0.19$~mm travelled per prediction window so that
turns are not averaged out by the windowing, and turn magnitudes drawn
with probability $0.5$ from $[120^\circ, 170^\circ]$ (near reversals)
and otherwise from $[30^\circ, 120^\circ]$; reflections at the
scan-area boundary confine the path and contribute additional corners.
The diagonal test geometry of Section~\ref{sec:results} is excluded from
both ensembles, so it remains geometrically unseen for both surrogates.
Figure~\ref{fig:path-profiles} shows some of the paths used for the training and testing.

\begin{figure}[tb]
    \centering
    \includegraphics[width=0.5\linewidth]{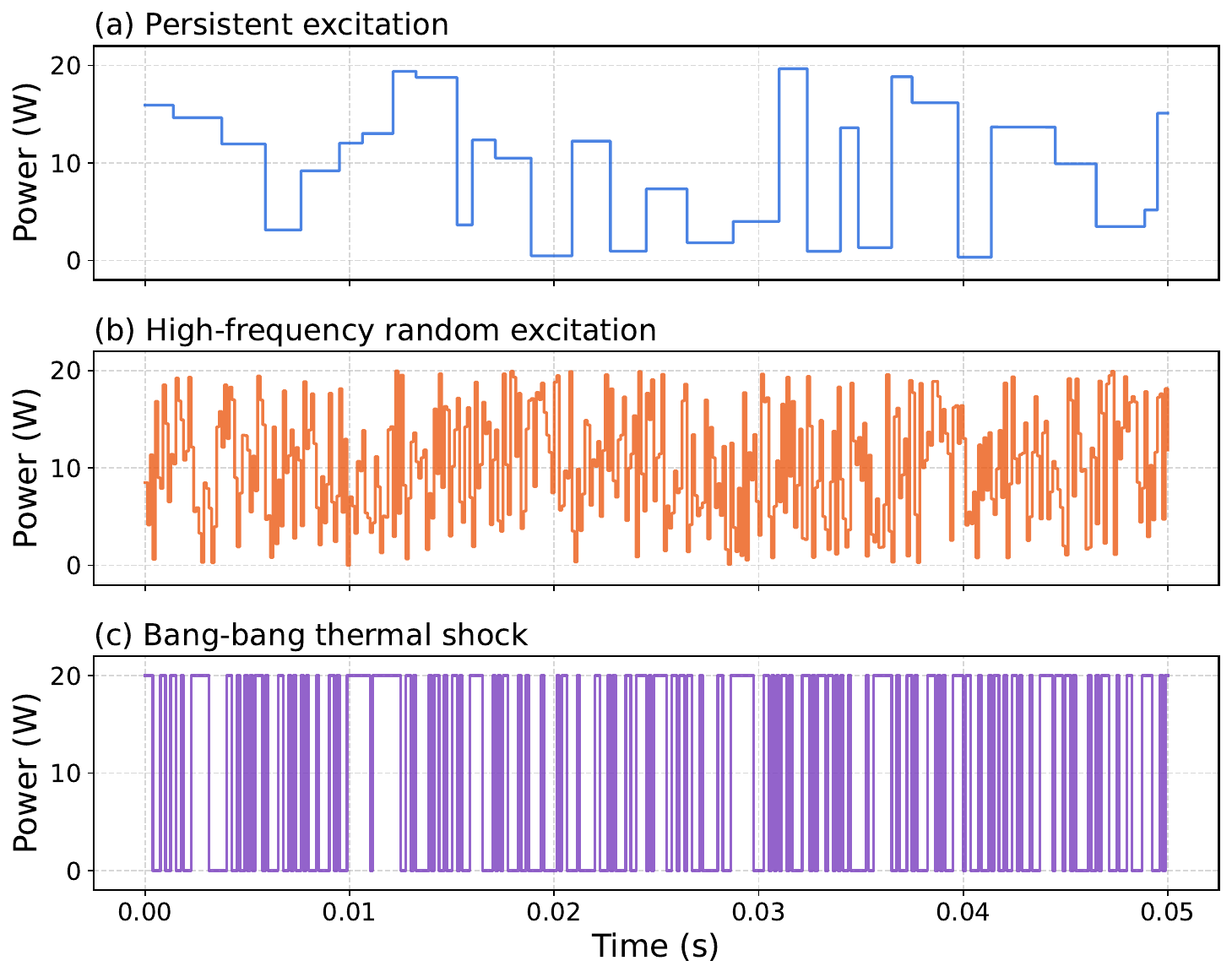}
    \caption{The three laser-power excitation classes used for
    training-data generation: (a)~persistent excitation, (b)~high-frequency
    random excitation, and (c)~bang-bang thermal shock. Amplitudes span
    $[0, P_{\max}] = [0, 20]$~W.}
    \label{fig:power-profiles}
\end{figure}

\paragraph{Stochastic power profiles}
Each run is driven by one of three excitation classes
(Fig.~\ref{fig:power-profiles}), sampled in equal shares and identically for
both ensembles: (i)~\emph{persistent excitation} -- piecewise-constant
power held for random durations of 5--20 steps with amplitudes drawn
from $\mathcal{U}(0, P_{\max})$, probing the slow thermal modes;
(ii)~\emph{high-frequency random excitation} -- an independent amplitude
per step, probing the fast response; and (iii)~\emph{bang-bang thermal
shock} -- random switching between $0$ and $P_{\max}$, probing extreme
rate transients. Mixing excitation classes follows standard
system-identification practice on persistently exciting
inputs~\cite{willems2005persistency}: the closed-loop controller issues rapid power adjustments, and a training set dominated by slowly varying inputs
would leave the high-frequency response poorly identified.

\begin{figure}[tb]
    \centering
    \begin{tikzpicture}[
        >={Latex[length=2mm]},
        font=\small,
        cell/.style={draw=cData!60!black, minimum width=8.5mm,
                     minimum height=8.5mm, inner sep=0pt},
        lbl/.style={font=\scriptsize, align=center}
    ]
        \foreach \i [count=\c from 0] in {1,2,3,4,5,6,7,8,9} {
            \node[cell, fill=cData!8] (s\i) at (\c*0.9,0) {\scriptsize $t_{\i}$};
        }
        \node at (9.3,0) {$\cdots$};
        \node[cell, fill=cData!8] (sN) at (10.0,0) {\scriptsize $t_{400}$};

        \draw[decorate, decoration={brace, amplitude=4pt}, thick, cInput]
            (s1.north west) -- (s1.north east)
            node[midway, above=5pt, lbl, text=cInput!70!black]
            {input state $\mathbf{y}_{1}$};
        \draw[decorate, decoration={brace, amplitude=4pt}, thick, cGate]
            (s2.north west) -- (s6.north east)
            node[midway, above=5pt, lbl, text=cGate!80!black]
            {target $\mathbf{s}_{1} = T_{\max}(t_{2{:}6})$};
        \foreach \i in {1} { \node[cell, fill=cInput!18] at (s\i) {\scriptsize $t_{\i}$}; }
        \foreach \i in {2,3,4,5,6} { \node[cell, fill=cGate!18] at (s\i) {\scriptsize $t_{\i}$}; }

        \draw[decorate, decoration={brace, mirror, amplitude=4pt}, thick, cModel]
            (s2.south west) -- (s2.south east)
            node[midway, below=5pt, lbl, text=cModel!75!black]
            {$\mathbf{y}_{2}$};
        \draw[decorate, decoration={brace, mirror, amplitude=4pt}, thick, cModel]
            (s3.south west) -- (s7.south east)
            node[midway, below=5pt, lbl, text=cModel!75!black]
            {$\mathbf{s}_{2} = T_{\max}(t_{3{:}7})$};

        \draw[->, thick, cOut!70!black]
            ([yshift=-9mm]s1.south) -- ([yshift=-9mm]s2.south)
            node[midway, below, lbl, text=cOut!70!black] {stride $=1$};
    \end{tikzpicture}
    \caption{Overlapping sliding-window extraction with unit stride. Each
    window pairs the input state $\mathbf{y}_k$ at $t_k$ (blue) with the
    branch actuation/kinematics over $[t_k, t_{k+\Nh}]$ and the target
    $\mathbf{s}_k = T_{\max}(t_{k+1:k+\Nh})$ (orange). Consecutive windows
    (purple) are shifted by one step, so a $400$-frame trajectory yields
    $400 - \Nh = 395$ supervised windows.}
    \label{fig:windowing}
\end{figure}
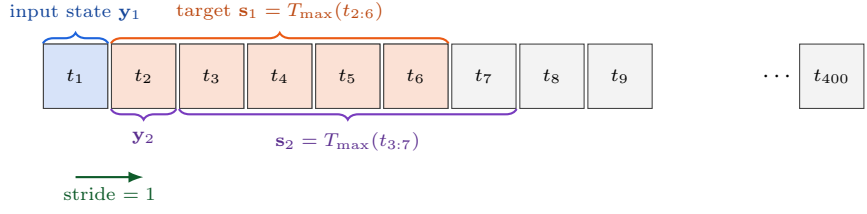

\paragraph{Windowing, scaling, and training}
Each trajectory is sliced into overlapping windows with unit stride (Fig.~\ref{fig:windowing}): a single $400$-step simulation produces $N_t - \Nh = 395$ supervised triples
$(\mathbf{u}_k, \mathbf{y}_k, \mathbf{s}_k)$ with target
$\mathbf{s}_k = \bigl( \Tmax(t_{k+1}), \dots, \Tmax(t_{k+\Nh}) \bigr)^\top$,
for a total of $\approx 1.26 \times 10^{5}$ samples. Inputs and targets are standardized per feature
(zero mean, unit variance), the data are shuffled, and an $80{:}20$
train/validation split is drawn with a fixed seed. 
Both surrogates share an identical architecture
(Section~\ref{sec:arch}) and training protocol: mean squared error
on normalized targets minimized with Adam~\cite{kingma2015adam} in
PyTorch~\cite{paszke2019pytorch} (learning rate $10^{-3}$,
weight decay $10^{-5}$, batch size $1024$), a plateau-based learning-rate
schedule (halving after $100$ stagnant epochs), and early stopping on the
validation loss with a $300$-epoch patience.
Any performance difference between the
two surrogates is therefore attributable solely to the geometry
composition of the training data.

\subsection{Offline surrogate validation}
\label{sec:offline-validation}

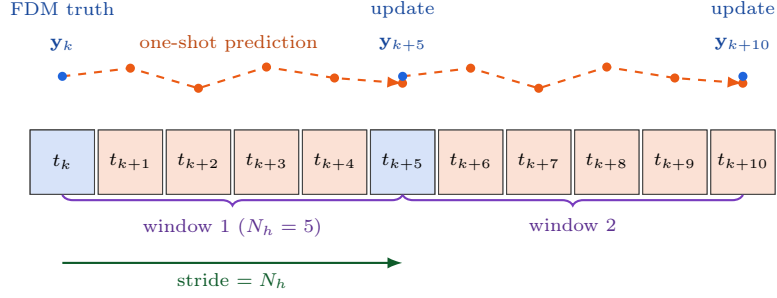
\begin{figure}[tb]
    \centering
    \begin{tikzpicture}[
        >={Latex[length=2mm]},
        font=\small,
        cell/.style={draw=cData!60!black, minimum width=8.5mm,
                     minimum height=8.5mm, inner sep=0pt},
        lbl/.style={font=\scriptsize, align=center}
    ]
        \foreach \i [count=\c from 0] in
            {k,{k{+}1},{k{+}2},{k{+}3},{k{+}4},{k{+}5},
             {k{+}6},{k{+}7},{k{+}8},{k{+}9},{k{+}10}} {
            \node[cell, fill=cData!8] (s\c) at (\c*0.9,0) {\scriptsize $t_{\i}$};
        }

        \foreach \i [count=\c from 0] in
            {k,{k{+}1},{k{+}2},{k{+}3},{k{+}4}} {
            \node[cell, fill=cGate!18] at (s\c) {\scriptsize $t_{\i}$};
        }
        \node[cell, fill=cInput!18] at (s0) {\scriptsize $t_{k}$};   
        \foreach \i [count=\c from 5] in
            {{k{+}5},{k{+}6},{k{+}7},{k{+}8},{k{+}9}} {
            \node[cell, fill=cGate!18] at (s\c) {\scriptsize $t_{\i}$};
        }
        \node[cell, fill=cInput!18] at (s5) {\scriptsize $t_{k{+}5}$}; 
        \node[cell, fill=cGate!18] at (s10) {\scriptsize $t_{k{+}10}$};

        \foreach \b in {0,5} {
            \draw[cGate, thick, dashed, ->]
                ([yshift=20pt]s\b.north)
                \foreach \i in {1,...,5} {
                    -- ([yshift={20 + 4.5*sin((\i)*137)},
                         xshift={(\i)*25.6}]s\b.north)
                };
            \foreach \i in {1,...,5} {
                \fill[cGate] ([yshift={20 + 4.5*sin((\i)*137)},
                              xshift={(\i)*25.6}]s\b.north) circle (0.06);
            }
        }
        \node[lbl, text=cGate!80!black, anchor=south]
            at ([yshift=9mm,xshift=-5mm]s3.north) {one-shot prediction};

        \foreach \c in {0,5,10} {
            \fill[cInput] ([yshift=20pt]s\c.north) circle (0.06);
        }
        \node[lbl, text=cInput!70!black, anchor=south]
            at ([yshift=26pt]s0.north)  {FDM truth\\ $\mathbf{y}_{k}$};
        \node[lbl, text=cInput!70!black, anchor=south]
            at ([yshift=26pt]s5.north)  {update\\ $\mathbf{y}_{k+5}$};
        \node[lbl, text=cInput!70!black, anchor=south]
            at ([yshift=26pt]s10.north) {update\\ $\mathbf{y}_{k+10}$};

        \draw[decorate, decoration={brace, mirror, amplitude=4pt}, thick, cModel]
            (s0.south) -- (s5.south)
            node[midway, below=5pt, lbl, text=cModel!75!black]
            {window 1 ($\Nh = 5$)};
        \draw[decorate, decoration={brace, mirror, amplitude=4pt}, thick, cModel]
            (s5.south) -- (s10.south)
            node[midway, below=5pt, lbl, text=cModel!75!black]
            {window 2};

        \draw[->, thick, cOut!70!black]
            ([yshift=-9mm]s0.south) -- ([yshift=-9mm]s5.south)
            node[midway, below, lbl, text=cOut!70!black] {stride $=\Nh$};
    \end{tikzpicture}
    \caption{Non-overlapping validation protocol. The trunk state is
    read from the FDM ground truth only at window boundaries (blue): within each window the surrogate predicts all $\Nh$ steps in one shot without access to the intermediate truth (orange).}
    \label{fig:protocol}
\end{figure}

Overlapping training windows share $\Nh - 1$ states with their neighbors, so
an evaluation under the same protocol is not adequate. We
therefore validate on strictly \emph{non-overlapping} windows
(Fig.~\ref{fig:protocol}): at step $k$ the trunk state $\mathbf{y}_k$ is
read from the FDM ground truth, the network predicts the next $\Nh$ values
of $\Tmax$ with no access to the intermediate truth, and only at step
$k + \Nh$ is the trunk state updated again with the FDM value. Prediction errors thus accumulate
freely within each window, exposing exactly the multi-step behavior on which
the MPC relies.

The protocol is applied to trajectories whose power profiles and, where
stated, geometries were never used in training: seen geometries under unseen
excitation separate representation accuracy from memorization, while the
diagonal path tests out-of-distribution generalization. Accuracy is reported
as the per-horizon-step and pooled root-mean-square error and the maximum
absolute error against the FDM; the corresponding results, covering all
three power-profile classes, are collected in Section~\ref{sec:results}.

\section{Surrogate-embedded model predictive control}
\label{sec:mpc}

\subsection{Differentiable NLP embedding}
\label{sec:embedding}

Problem~\eqref{eq:mpc} is solved with the interior-point solver
IPOPT through the CasADi framework~\cite{andersson2019casadi}, which
requires the dynamics constraint~\eqref{eq:mpc-dyn} to be available as a
smooth symbolic expression. To this end, the trained PyTorch weights of the
MS-DeepONet are extracted and the network~\eqref{eq:msdeeponet} is rebuilt
layer by layer as a CasADi computational graph: the per-feature
standardization of Section~\ref{sec:data}, both subnetworks, the
branch--trunk inner product, and the inverse output scaling are composed
into a single symbolic function
$\mathcal{F}_{\text{surr}} : \mathbb{R}^{5\Nh} \times \mathbb{R}^{1+\Nh}
\to \mathbb{R}^{\Nh}$ whose exact first- and second-order derivatives are
obtained by automatic differentiation. The optimizer therefore
differentiates the identical map that was trained, rather than a
finite-difference approximation of it.

One modification is required. The ReLU activation
$\sigma(z) = \max(0, z)$ is continuous but not differentiable at $z = 0$,
and its second derivative vanishes almost everywhere; both properties
degrade the Hessian information on which interior-point methods rely and
were observed to stall or destabilize IPOPT. During embedding, and only there, each ReLU is replaced by the smooth algebraic approximation
\begin{equation}
    \tilde{\sigma}(z)
    \;=\;
    \tfrac{1}{2}\Bigl( z + \sqrt{z^{2} + \epsilon} \Bigr),
    \qquad \epsilon = 10^{-6},
    \label{eq:smooth-relu}
\end{equation}
which is $C^{\infty}$, satisfies
$\tilde{\sigma}(z) \to \sigma(z)$ as $\epsilon \to 0$, and deviates from the
exact ReLU by at most $\sqrt{\epsilon}/2 = 5 \times 10^{-4}$ (attained at
$z = 0$). The trained weights are used unchanged; the perturbation of the
network output is negligible relative to the surrogate's prediction error,
while the resulting $C^{2}$ constraint function restores reliable Newton
steps in the solver.

Because the optimizer acts on the surrogate's derivatives rather than its
values alone, we assess the \emph{gradient fidelity} of the embedded network
on a representative seen path (spiral). At operating points sampled along
the path we compare the analytic sensitivity
$\partial \widehat{T}_{\max,k+1} / \partial P_{k+1}$ of the CasADi graph,
obtained by automatic differentiation, against a central finite difference
of the FDM, which serves as the true physical sensitivity. Two observations
follow. First, the surrogate reproduces the \emph{sign} of the FDM
sensitivity in $100\%$ of the tested points, so the descent directions
supplied to IPOPT are always physically consistent. Second, the magnitude is
captured only approximately: with a typical sensitivity of
$\approx 7$~K\,W$^{-1}$, the mean absolute discrepancy is $2.8$~K\,W$^{-1}$
(median $2.5$~K\,W$^{-1}$). For the
interior-point solver this suffices, since its line search rescales the step
length and thus tolerates moderate errors in gradient magnitude provided the
direction is correct.

\subsection{Closed-loop implementation}
\label{sec:closedloop-impl}
 
The receding-horizon loop realizes the architecture of
Fig.~\ref{fig:closedloop} as summarized in
Algorithm~\ref{alg:mpc}. At step $k$, the top-surface temperature field
of the FDM, standing in for a thermal-camera image, is queried at the
current peak and at the $N_h$ known future path locations to assemble
the trunk state $\mathbf{y}_k$ of Eq.~\eqref{eq:trunk}. No quantity
beyond this surface information is fed back, so the controller relies
only on measurements available in practice. The NLP~\eqref{eq:mpc} is
then solved, the first optimal power $P^{\star}_{k+1}$ is applied, the
FDM advances the true state by one step $\Delta t$, and the loop repeats.

\begin{algorithm}[tb]
\caption{Surrogate-embedded MPC}
\label{alg:mpc}
\begin{algorithmic}[1]
\Require trained MS-DeepONet weights; scan path $\bm{\xi}(t)$; horizon $\Nh$;
         bounds $[T_{\min}, T_{\max}^{\lim}]$, $P_{\max}$
\State build symbolic $\mathcal{F}_{\text{surr}}$ with smooth
       ReLU~\eqref{eq:smooth-relu}; instantiate IPOPT on
       problem~\eqref{eq:mpc}
\State initialize plant at $T(\mathbf{x}, 0) = T_0$; set $P_0^{\star}$ to a
       nominal power
\For{$k = 0, 1, 2, \dots$}
    \State read surface field of the FDM plant (emulated thermal camera)
    \State assemble $\mathbf{y}_k$: $\Tmax(t_k)$ and
           $T(\bm{\xi}_{k+i}, t_k)$, $i = 1{:}\Nh$
           \Comment{Eq.~\eqref{eq:trunk}}
    \State assemble path kinematics of $\mathbf{u}_k$ from $\bm{\xi}$
           \Comment{Eq.~\eqref{eq:branch}}
    \State warm-start with $P_k^{\star}$ and slack guess
           $\max(0,\, T_{\min} - \Tmax(t_k))$
    \State solve NLP~\eqref{eq:mpc} $\rightarrow$
           $\mathbf{P}_k^{\star}, \bm{\epsilon}_k^{\star}$
    \State apply first element $P_{k+1}^{\star}$; advance FDM plant by
           $\Delta t$
\EndFor
\end{algorithmic}
\end{algorithm}
 
Three implementation details matter for solver robustness. First, the
decision variables are normalized: powers as $\bar{P} = P / P_{\max}$
and slacks as $\bar{\epsilon} = \epsilon / \epsilon_{\max}$ with
$\epsilon_{\max} = 500$~K, so that all variables are of order one and
the $\ell_1$ penalty weight retains its intended dominance. Second, the
solver is warm-started: the previous optimal power (held constant across
the horizon) initializes $\mathbf{P}_k$, and the slack guess is set to
the currently observed lower-bound violation
$\max\big(0,\, T_{\min} - T_{\max}(t_k)\big)$, which is nonzero only
during the initial heat-up transient. Third, IPOPT is run with its
limited-memory (L-BFGS) quasi-Newton Hessian approximation and an
iteration cap of $500$; exact Hessians of the composed network graph offer
no observed benefit at this scale while increasing the per-iteration
cost. With these settings the solver converges at every step of all
closed-loop experiments reported in Section~\ref{sec:results}.

\paragraph{Constraint tightening against surrogate error}
The hard upper bound~\eqref{eq:mpc-hard} is enforced on the surrogate's
prediction $\widehat{T}_{\max}$, not on the true plant temperature.
Consequently, whenever the surrogate under-predicts, the physical
temperature can transiently exceed $T^{\lim}_{\max}$ even though the
optimizer reports a feasible plan. 
Section~\ref{sec:results-cl-ood}
shows that the corner-augmented surrogate, despite eliminating the
large corner-localized failure of the baseline model, retains a
residual in-distribution prediction error that produces such small
transient exceedances on the "smooth" seen paths. We absorb this residual
error by a calibrated safety margin: the upper bound supplied to the
optimizer is tightened to
\begin{equation}
    \widehat{T}_{\max,k+i} \;\le\; T^{\lim}_{\max} - \delta,
    \qquad i = 1, \dots, N_h.
    \label{eq:tightening}
\end{equation}
Because only under-prediction can cause a violation of the upper bound
(an over-predicting surrogate errs in the safe direction), the margin
is calibrated one-sidedly: $\delta$ is set to the 95th percentile of
the one-step under-prediction residual
$\max\!\big(0,\, T^{\mathrm{FDM}}_{\max} - \widehat{T}_{\max}\big)$ of
the deployed corner-augmented surrogate, evaluated on a held-out $20\%$
validation split of the windowed dataset ($25{,}280$ windows, drawn
with a fixed seed and pooling all three excitation classes). This
yields $13.2$~K, which we round down to $\delta = 13$~K; the
corresponding symmetric margin (95th percentile of the absolute
residual, $18.5$~K) would consume nearly half of the $40$~K process
window and is unnecessarily conservative for a one-sided constraint.
The calibration uses offline data only, no plant information beyond
that already used for training, and adds no computational cost to the
online problem. Two limitations are stated explicitly: the margin is
calibrated on the one-step \emph{in-distribution} residual, so it is
not designed to cover the larger extrapolation error that a
geometrically out-of-distribution feature may still induce; and the
quantile-based construction carries no formal robustness guarantee. A
rigorous tube-based or distributionally robust treatment of the
learned-surrogate error is left to future work. The closed-loop effect
of the margin, including its conservatism cost, is quantified in
Section~\ref{sec:results-cl}.

\section{Results}
\label{sec:results}

\subsection{Offline prediction accuracy}
\label{sec:offline-results}
 
\begin{figure}[tb]
    \centering
    \includegraphics[width=0.9\linewidth]{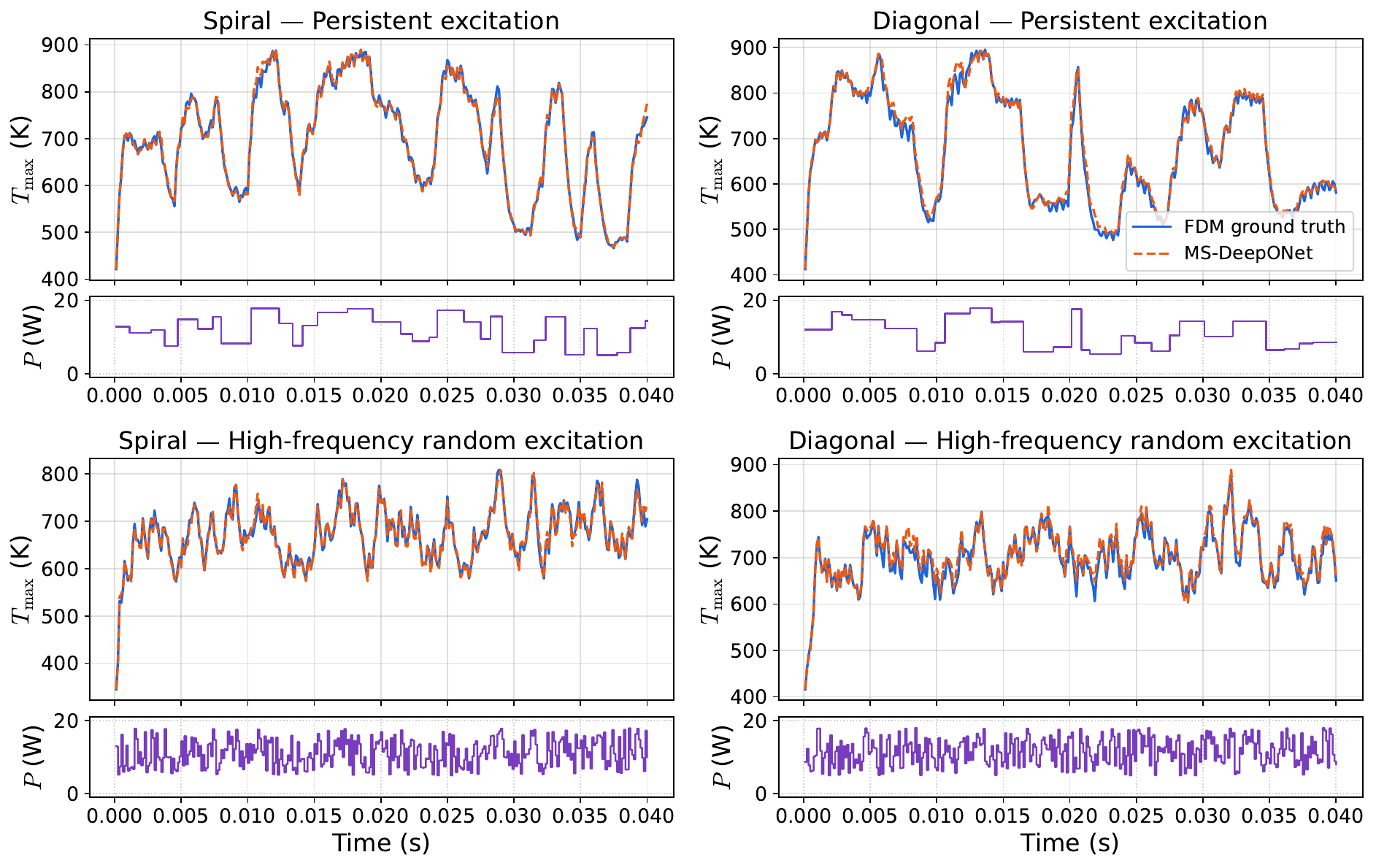}
    \caption{Offline multi-step validation of the
    corner-augmented MS-DeepONet: multi-step prediction vs. FDM on the seen (spiral) and unseen (diagonal) geometry.}
    \label{fig:offline}
\end{figure}

We evaluate both surrogates under the non-overlapping blind protocol of Section~\ref{sec:offline-validation} on the seen geometries (vertical,
spiral) and on the geometrically unseen diagonal path, driven by newly
sampled realizations of all three excitation classes of
Section~\ref{sec:data} (persistent, high-frequency random, and
bang-bang). To make the comparison exact, the two surrogates are scored
on identical FDM trajectories: for each (geometry, excitation)
pair the plant is simulated once and both models predict the same
ground truth, so every difference in Table~\ref{tab:offline} is
attributable to the model alone. Each entry is averaged over three
independent excitation realizations; the run-to-run standard deviation
of the horizon RMSE is at most $1$~K, small compared to all reported
differences. Figure~\ref{fig:offline} illustrates the
resulting $\Tmax$ tracking for representative cases.
 
\begin{table}[bt]
\centering
\caption{Offline multi-step prediction accuracy. Horizon
RMSE per excitation class, its mean over the three classes, and
the worst absolute error across all runs, in K.}
\label{tab:offline}
\scriptsize
\begin{tabular}{llccccc}
\toprule
& & \multicolumn{3}{c}{Horizon RMSE (K)} & & \\
\cmidrule(lr){3-5}
Scan path & Surrogate & Persistent & HF-random & Bang-bang & Mean & $\max|e|$ \\
\midrule
Vertical (seen)   & Baseline    & 2.6  & 3.4  & 5.3  & 3.7  & 27.7 \\
                  & Corner-aug. & 4.7  & 5.4  & 7.7  & 5.9  & 37.2 \\
\midrule
Spiral (seen)     & Baseline    & 3.8  & 4.1  & 7.0  & 4.9  & 34.9 \\
                  & Corner-aug. & 8.0  & 8.2  & 10.1 & 8.8  & 41.9 \\
\midrule
Diagonal (unseen) & Baseline    & 12.1 & 13.6 & 13.4 & 13.0 & 79.9 \\
                  & Corner-aug. & 13.7 & 14.2 & 14.6 & 14.2 & 71.2 \\
\bottomrule
\end{tabular}
\end{table}

Three observations structure the comparison. First, on the seen
geometries the baseline surrogate is accurate: its pooled horizon RMSE
is $3.7$~K on the vertical path and $4.9$~K on the spiral, i.e.\ well
below $1\%$ of the regulated temperature level. The error increases
monotonically with the frequency content of the excitation (persistent
$<$ high-frequency random $<$ bang-bang), as the fast thermal
transients are the hardest to resolve, and is slightly larger on the
spiral, whose continuously rotating velocity vector exercises a broader
range of source kinematics than the piecewise-straight vertical path.
The per-horizon-step errors are essentially flat from the first to the
fifth prediction step (e.g.\ $2.7$ to $3.0$~K on the vertical path
under persistent excitation), confirming that the one-shot formulation
does not accumulate error within the blind window; the contrast with
recursive predictors is quantified in Section~\ref{sec:ablations}.
 
Second, the corner augmentation degrades the in-distribution accuracy,
quantifying the trade-off announced in Section~\ref{sec:data}: the
pooled RMSE of the corner-augmented surrogate rises to $5.9$~K on the
vertical path ($+59\%$) and $8.8$~K on the spiral ($+78\%$), with worst
errors growing from $28$ to $37$~K and from $35$ to $42$~K,
respectively. Rebalancing a fixed-size training budget toward the
corner regime dilutes the density of smooth-path samples, and the
surrogate cedes part of its in-distribution accuracy in exchange for
coverage of the reversal kinematics.
 
Third, and in contrast to the seen-path cost, the aggregate open-loop benefit of the augmentation on the unseen diagonal is barely
visible: pooled over excitation classes, the two surrogates are
statistically indistinguishable ($13.0$ vs.\ $14.2$~K), and the
worst-case error is reduced only modestly ($79.9 \to 71.2$~K). Averaged
over broadband stochastic excitation, both models appear equally
(in)adequate on the unseen geometry. This is not a deficiency of the
augmentation but a limitation of trajectory-averaged open-loop
validation. The failure mode that matters for control is
regime-specific: it arises when the source re-enters strongly
pre-heated material at the sharp reversal \emph{while the temperature
is held near the upper operating bound} -- a state that a randomly
excited trajectory visits only rarely and transiently, so its signature
is diluted in aggregate metrics, but that the closed-loop controller,
which regulates $\Tmax$ into the band $[760, 800]$~K by design,
occupies permanently. Under closed-loop conditions the two surrogates
behave drastically differently at exactly this regime
(Section~\ref{sec:results-cl-ood}): the baseline under-predicts the
corner accumulation by up to $91$~K and causes a $51$~K violation of
the hard bound, whereas the corner-augmented model contains the
worst-case exceedance to $1.4$~K. Open-loop accuracy is therefore a
necessary but not a sufficient indicator of control-readiness, and the
closed-loop evaluation of Sections~\ref{sec:results-cl-seen}--\ref{sec:results-cl}
is the decisive test.

\subsection{Closed-loop MPC on seen scan paths}
\label{sec:results-cl-seen}
 
We first deploy the baseline surrogate (trained on
$\mathcal{D}_{\mathrm{base}}$) in the receding-horizon loop of
Algorithm~\ref{alg:mpc} on the seen geometries. All closed-loop
experiments span $N_{\mathrm{sim}} = 320$ control steps ($0.04$~s), and
all statistics reported below are evaluated against the FDM ground
truth over the settled window $k \ge 10$ ($t \ge 1.25$~ms), which
excludes the initial heat-up ramp during which the plate is heated from
the ambient temperature toward the operating band and the soft lower
bound is necessarily violated. Table~\ref{tab:closedloop} collects the
constraint-satisfaction, prediction-accuracy, and control-effort
metrics of all closed-loop experiments; the spiral results, which are
qualitatively identical to the vertical ones, are deferred to
\ref{app:spiral}.

\begin{table}[tb]
\centering
\caption{Closed-loop statistics of all MPC experiments
($N_{\mathrm{sim}} = 320$ steps), evaluated against the FDM ground
truth over the settled window $k \ge 10$. Overshoot is the largest
exceedance of the true hard bound $T^{\lim}_{\max} = 800$~K;
$n_{>800}$ counts the control steps in violation. The residual columns
report the one-step prediction error
$\widehat{T}_{\max} - T^{\mathrm{FDM}}_{\max}$; $\bar{e}$ is its signed
mean (positive = over-prediction). $\overline{|\Delta P|}$ is the mean
step-to-step power change. In the bottom group the optimizer enforces
the tightened bound $T^{\lim}_{\max} - \delta = 787$~K; scoring remains against $800$~K.}
\label{tab:closedloop}
\scriptsize
\begin{tabular}{llcccccccc}
\toprule
& & \multicolumn{3}{c}{Hard bound ($>800$~K)} & \multicolumn{3}{c}{Residual (K)} & \multicolumn{2}{c}{Power (W)} \\
\cmidrule(lr){3-5} \cmidrule(lr){6-8} \cmidrule(lr){9-10}
Surrogate & Scan path & Overshoot (K) & $n_{>800}$ & Dur.\ (ms) & RMS & $\max |e|$ & $\bar{e}$ & mean & $\overline{|\Delta P|}$ \\
\midrule
Baseline    & Vertical (seen)    & 0.0  & 0  & 0.00 & 1.8  & 11.1 & $-0.7$ & 13.85 & 0.23 \\
Baseline    & Spiral (seen)      & 0.0  & 0  & 0.00 & 3.6  & 14.0 & $+0.0$ & 14.47 & 0.40 \\
Baseline    & Diagonal (unseen)  & 51.1 & 22 & 2.75 & 17.2 & 91.1 & $+3.3$ & 13.94 & 0.65 \\
\midrule
Corner-aug. & Vertical (seen)    & 5.8  & 1  & 0.13 & 3.6  & 17.1 & $+0.3$ & 13.82 & 0.40 \\
Corner-aug. & Spiral (seen)      & 2.7  & 5  & 0.63 & 6.9  & 23.9 & $+1.0$ & 14.40 & 0.54 \\
Corner-aug. & Diagonal (unseen)  & 1.4  & 1  & 0.13 & 14.8 & 37.7 & $+9.3$ & 13.64 & 0.58 \\
\midrule
Corner-aug.\ $+\,\delta$ & Vertical (seen)   & 0.0 & 0 & 0.00 & 3.6  & 17.1 & $+0.3$ & 13.82 & 0.41 \\
Corner-aug.\ $+\,\delta$ & Spiral (seen)     & 0.0 & 0 & 0.00 & 6.9  & 23.8 & $+1.0$ & 14.36 & 0.63 \\
Corner-aug.\ $+\,\delta$ & Diagonal (unseen) & 0.0 & 0 & 0.00 & 14.7 & 37.5 & $+9.1$ & 13.58 & 0.79 \\
\bottomrule
\end{tabular}
\end{table}
 
\begin{figure}[tb]
    \centering
    \includegraphics[width=0.95\linewidth]{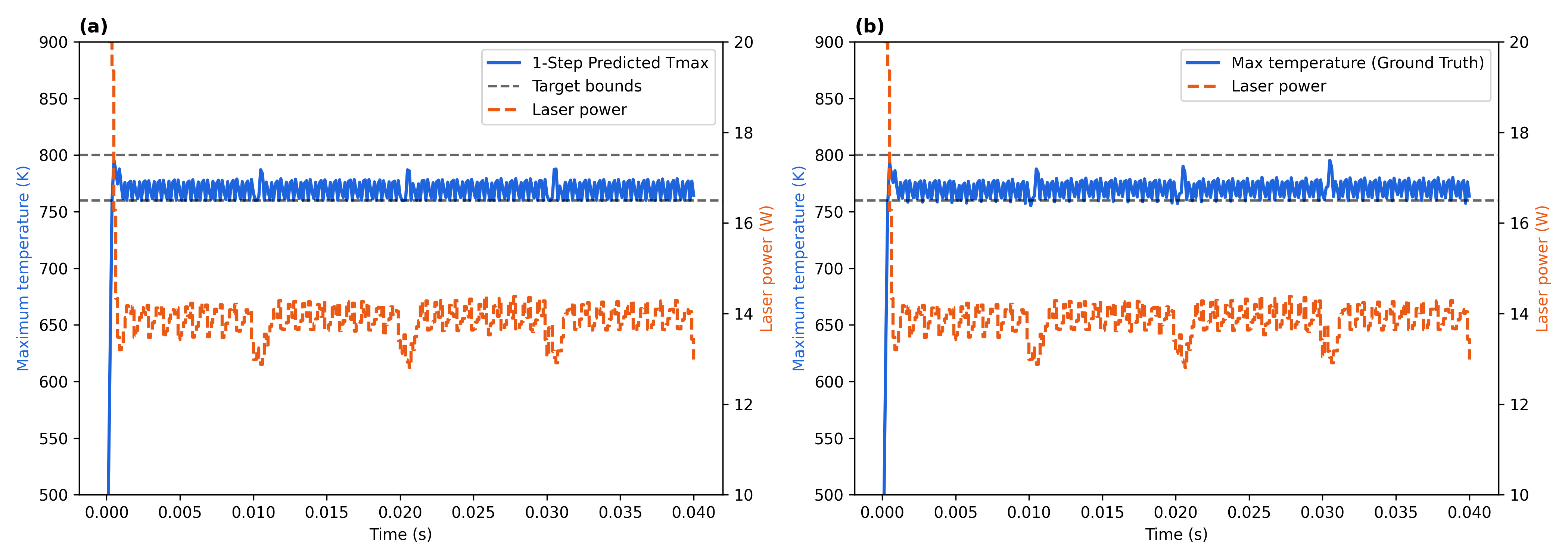}
    \caption{Closed-loop MPC on the seen vertical scan path with the
    baseline surrogate: (a)~surrogate one-step prediction of $\Tmax$
    and optimal laser power; (b)~FDM ground truth under the applied
    power sequence.}
    \label{fig:cl-base-vertical}
\end{figure}

Figure~\ref{fig:cl-base-vertical} shows the closed loop on the vertical
scan path. The controller settles within a few control steps and
regulates $\Tmax$ inside the prescribed band for the remainder of the
scan: the ground-truth peak temperature never exceeds the hard upper
bound (maximum $795.4$~K against the $800$~K limit), and excursions
below the soft lower bound are brief and shallow ($\le 4.8$~K),
consistent with the exact-penalty slack formulation, which tolerates
transient dips rather than forcing aggressive reheating. The
closed-loop one-step prediction residual quantifies how faithfully the
surrogate represents the plant along the trajectory that the controller
actually visits: on the vertical path the mean absolute residual is
$1.2$~K and the root-mean-square residual $1.8$~K, roughly $4\%$ of the
$40$~K width of the admissible band, with a worst instantaneous error
of $11.1$~K localized at the right-angle turn of the path. The
resulting power commands are smooth, with a mean step-to-step change of
$0.23$~W ($1.1\%$ of $P_{\max}$) and a maximum of $0.86$~W, reflecting
the rate penalty $w_{\Delta P}$ in the objective. On the spiral path
the picture is unchanged (zero hard-bound violations, RMS residual
$3.6$~K; \ref{app:spiral}); the moderately larger residual mirrors the
offline trend of Section~\ref{sec:offline-results} and is attributable
to the continuously rotating velocity vector of the spiral. On the
geometries it was trained on, the baseline surrogate is thus an
adequate internal model for constraint-enforcing control: the
prediction error is an order of magnitude smaller than the band width,
and the hard bound is respected throughout.
 
\subsection{Out-of-distribution analysis: unseen sharp-corner path}
\label{sec:results-cl-ood}
 
\begin{figure}[tb]
    \centering
    \includegraphics[width=0.95\linewidth]{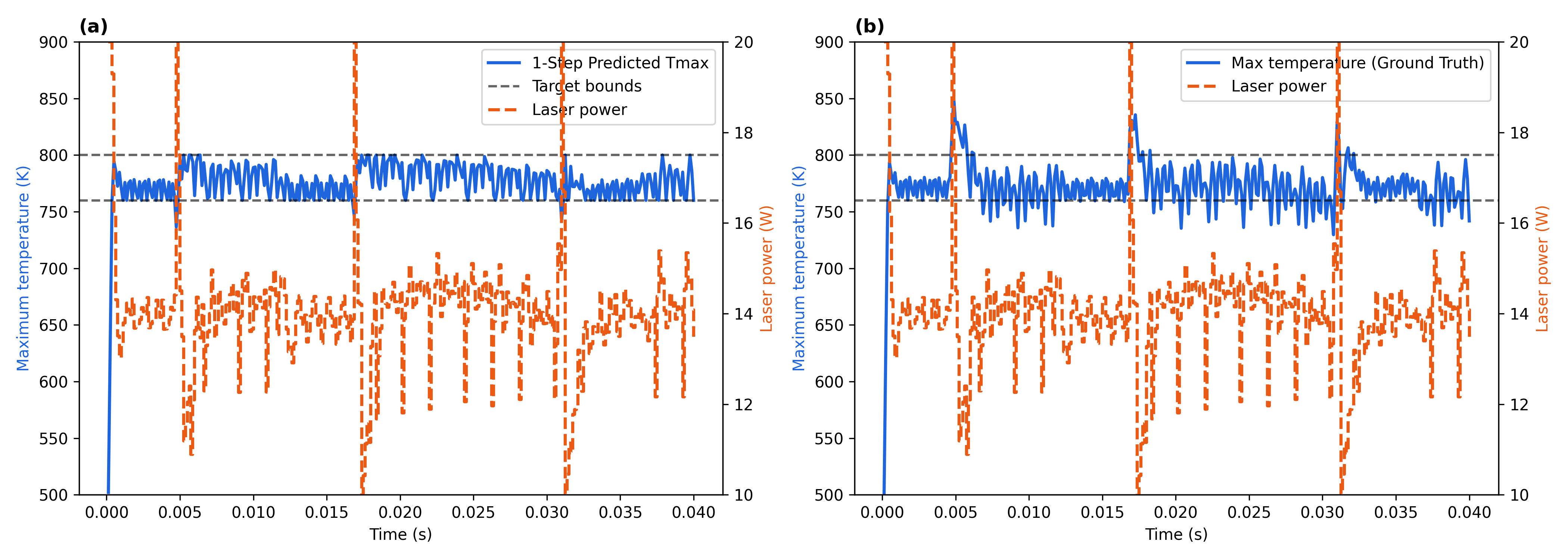}
    \caption{Closed-loop MPC on the unseen diagonal scan path with the
    baseline surrogate: (a)~surrogate prediction and optimal power;
    (b)~FDM ground truth. At the sharp direction reversal
    ($t \approx 4.9$~ms) the surrogate under-predicts the heat
    accumulation by up to $91$~K and the true temperature exceeds the
    hard bound by $51$~K.}
    \label{fig:cl-base-diagonal}
\end{figure}

\begin{figure}[tb!]
    \centering
    \begin{minipage}{0.49\linewidth}
        \centering
        \includegraphics[width=\linewidth]{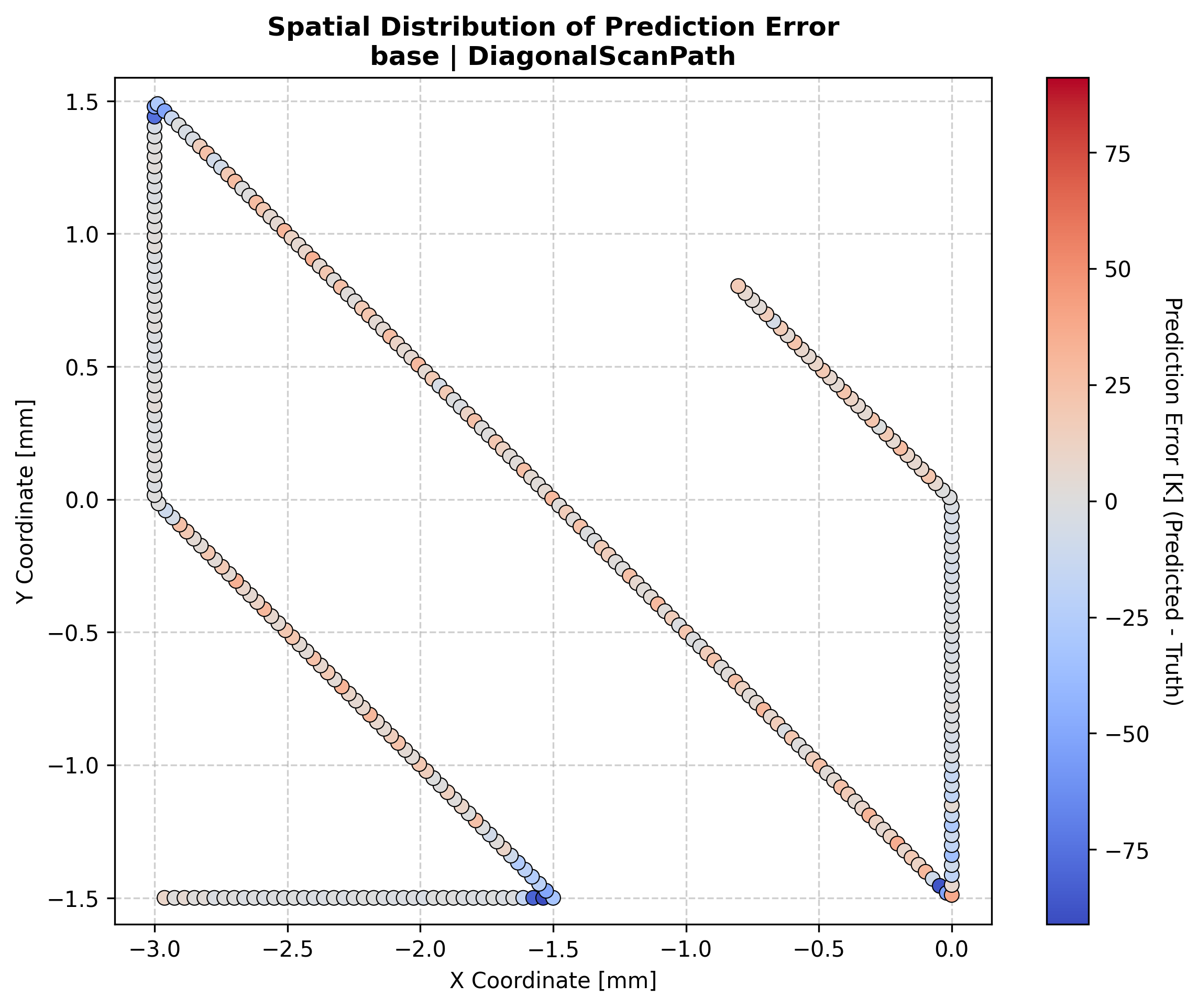}
        \\ (a)
    \end{minipage}
    \hfill
    \begin{minipage}{0.49\linewidth}
        \centering
        \includegraphics[width=\linewidth]{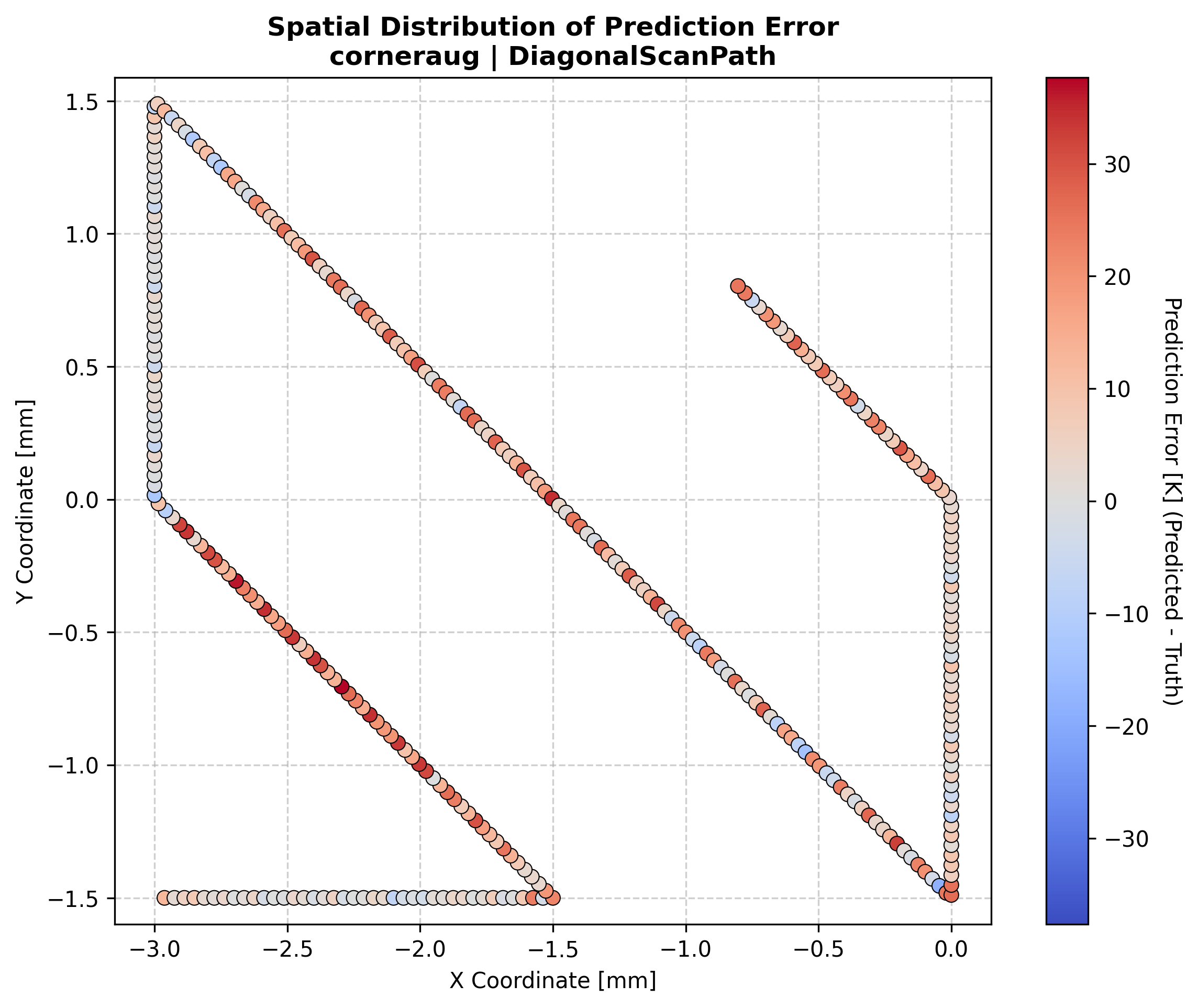}
        \\ (b)
    \end{minipage}
    \caption{Spatial distribution of the closed-loop one-step
    prediction residual along the unseen diagonal scan path:
    (a)~baseline surrogate -- the error is small along the straight
    segments but reaches $-91$~K at the sharp reversal vertex (note the
    color scale); (b)~corner-augmented surrogate -- the
    corner-localized under-prediction is eliminated and the residual
    is bounded by $\pm 38$~K with a predominantly positive
    (conservative) bias.}
    \label{fig:cl-spatial-diagonal}
\end{figure}

\paragraph{Failure of the baseline surrogate}
Deploying the same baseline controller on the geometrically unseen
diagonal scan path exposes the out-of-distribution failure mode
anticipated in Section~\ref{sec:data}. As shown in
Fig.~\ref{fig:cl-base-diagonal}, the controller tracks the band
adequately along the straight diagonal segments, but at the sharp
direction reversal of the path the ground-truth temperature overshoots
the hard bound by $51.1$~K, reaching $851.1$~K, with the violation
persisting for $22$ control steps ($2.75$~ms). The mechanism is read
off directly from the one-step residual: at the reversal vertex,
located at $(-1.5, -1.5)$~mm and reached at $t \approx 4.9$~ms, the
surrogate under-predicts the peak temperature by up to $91.1$~K. As the
source turns back into its own thermal wake, the true heat accumulation
exceeds anything represented in the smooth-path training distribution;
the surrogate, extrapolating from that distribution, reports a
feasible-looking plan, and the optimizer consequently fails to cut the
power in time. The spatial distribution of the residual
(Fig.~\ref{fig:cl-spatial-diagonal}a) confirms that the failure is
strictly localized: the error is small and unbiased along the straight
segments and collapses to a large negative value precisely at the
corner. Averaged over the whole trajectory the residual statistics
(RMS $17.2$~K) understate the problem; it is the localized, signed
error at the corner that breaks the constraint, which is why we report
the maximum error and its location alongside the aggregate metrics in
Table~\ref{tab:closedloop}.
 
\begin{figure}[tb]
    \centering
    \includegraphics[width=0.95\linewidth]{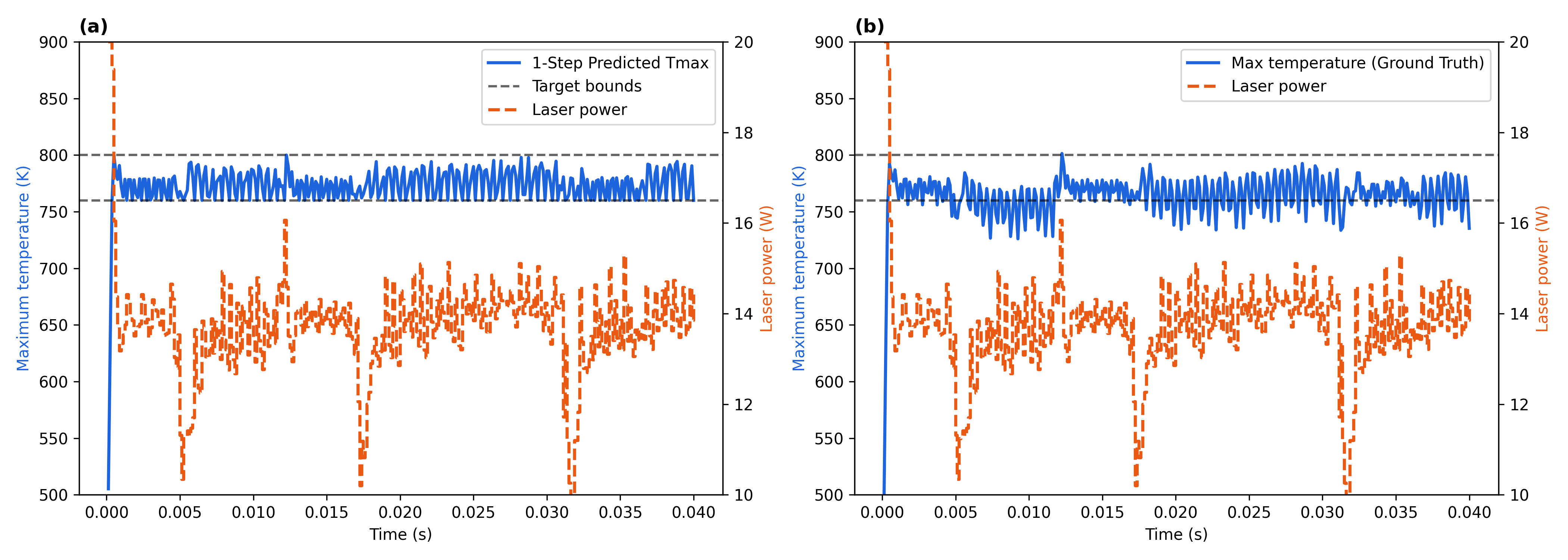}
    \caption{Closed-loop MPC on the unseen diagonal scan path with the
    corner-augmented surrogate (controller unchanged). The worst-case
    exceedance of the hard bound is reduced from $51.1$~K to $1.4$~K.}
    \label{fig:cl-corneraug-diagonal}
\end{figure}

\paragraph{Effect of corner-rich training data}
Replacing the surrogate by its corner-augmented counterpart (trained on
$\mathcal{D}_{\mathrm{corner}}$), with the controller and all weights
unchanged, removes the corner failure. On the same diagonal path
(Fig.~\ref{fig:cl-corneraug-diagonal}) the maximum ground-truth
temperature drops from $851.1$~K to $801.4$~K, a single-step exceedance
of $1.4$~K, i.e.\ a $97\%$ reduction of the worst-case violation; the
worst residual shrinks from $-91.1$~K to $+37.7$~K and, importantly,
changes sign. The spatial error map
(Fig.~\ref{fig:cl-spatial-diagonal}b) shows that the corner-localized
under-prediction has been eliminated and the remaining residual is
distributed along the path with a predominantly positive
(over-predicting) bias, $\bar{e} = +9.3$~K on this path. For the hard
upper bound this bias acts in the conservative direction: an
over-predicting surrogate causes the controller to withdraw power
early rather than late. Its price is visible at the lower end of the
band: the mean applied power drops by $0.3$~W relative to the baseline
run, and the fraction of settled steps below the soft lower bound
roughly doubles (from $19\%$ to $36\%$), with excursions of up to
$34$~K. The augmentation thus converts a safety-critical hard-bound
violation into a benign, soft-bound conservatism, which is the
preferable direction for a thermally limited process.

\begin{figure}[tb]
    \centering
    \includegraphics[width=0.95\linewidth]{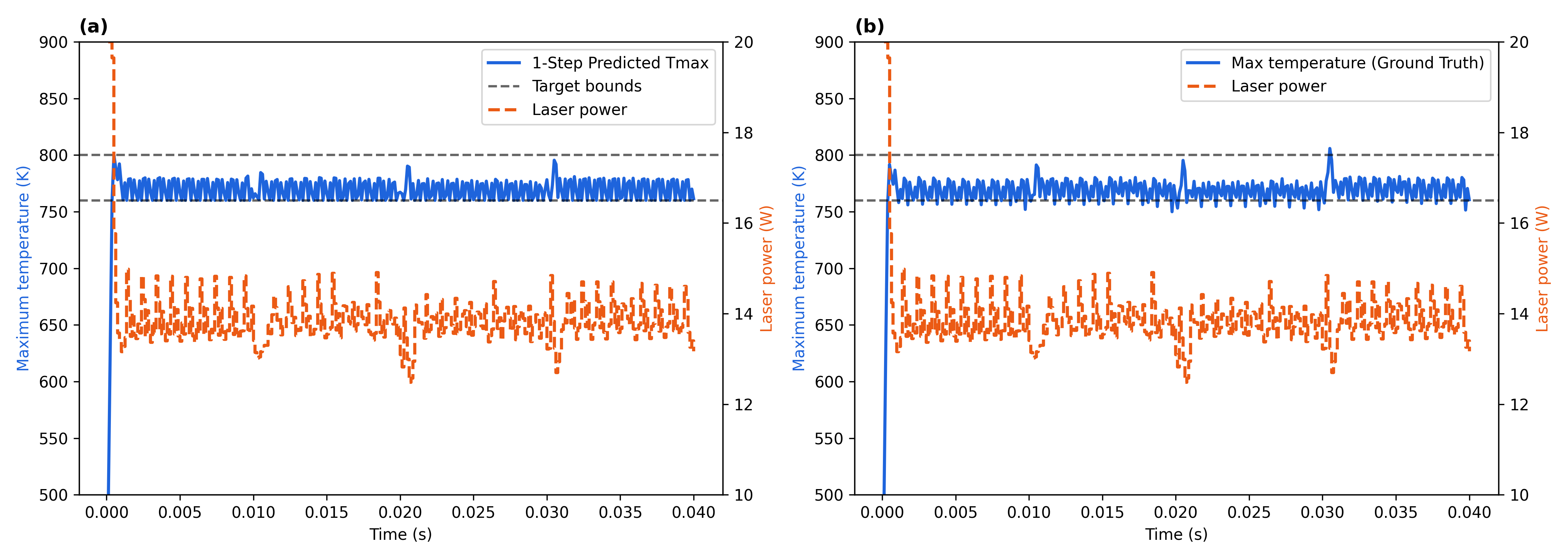}
    \caption{Closed-loop MPC on the seen vertical scan path with the
    corner-augmented surrogate. The increased in-distribution residual
    produces an isolated single-step exceedance of the hard bound by
    $5.8$~K (at $t \approx 30.5$~ms), motivating the calibrated
    constraint tightening of Eq.~\eqref{eq:tightening}.}
    \label{fig:cl-corneraug-vertical}
\end{figure}

\paragraph{Cost of the augmentation on seen scan paths}
The robustness gained at the corner is not free. Returning the
corner-augmented surrogate to the seen geometries quantifies the
trade-off announced in Section~\ref{sec:data}: on the vertical path the
closed-loop RMS residual doubles relative to the baseline surrogate
(from $1.8$~K to $3.6$~K) and the worst error grows from $11.1$~K to
$17.1$~K, concentrated at the same right-angle turn of the path; the
power commands become correspondingly less smooth
($\overline{|\Delta P|}$: $0.23 \to 0.40$~W). Figure~\ref{fig:cl-corneraug-vertical}
shows the consequence for constraint satisfaction: the ground truth now
exceeds the hard bound in isolated single-step events, by $5.8$~K on the
vertical path (one step, $0.125$~ms) and by up to $2.7$~K on the spiral
(five steps in total; \ref{app:spiral}). Rebalancing a fixed-size
training budget toward the previously under-represented corner regime
dilutes the density of smooth-path samples, and the surrogate cedes a
part of its in-distribution accuracy in exchange for coverage of the
reversal kinematics. The augmentation therefore changes the character
of the constraint-violation risk rather than eliminating it: the large,
localized out-of-distribution failure ($51$~K over $22$ steps) is
replaced by small, sporadic in-distribution exceedances ($\le 6$~K over
single steps) whose magnitude is proportional to the surrogate's
residual prediction error. Errors of this type are exactly what the
calibrated constraint tightening of Section~\ref{sec:closedloop-impl}
is designed to absorb, which we quantify next.
 
\subsection{Closed-loop control with constraint tightening}
\label{sec:results-cl}
 
\begin{figure}[tb]
    \centering
    \includegraphics[width=0.95\linewidth]{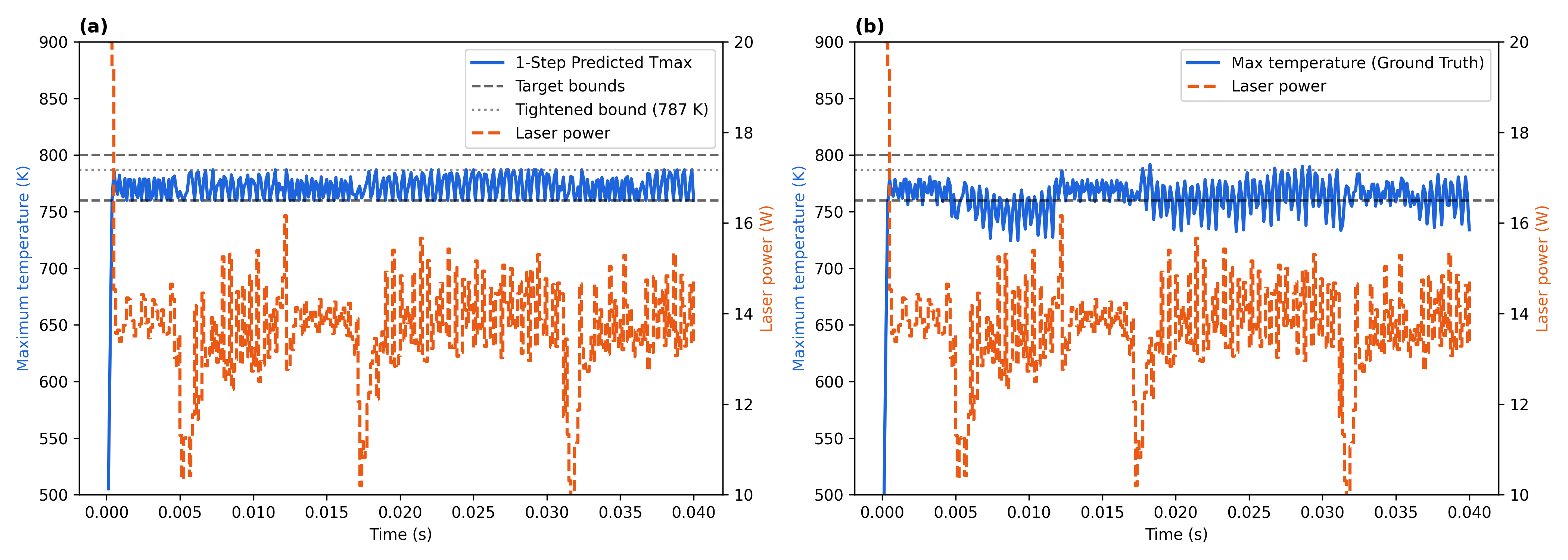}
    \caption{Closed-loop MPC on the unseen diagonal scan path with the
    corner-augmented surrogate and the tightened upper bound
    ($787$~K, dotted line): (a)~surrogate prediction and optimal power;
    (b)~FDM ground truth. The true peak temperature remains below the
    physical limit of $800$~K throughout ($\max\, 791.8$~K).}
    \label{fig:cl-tight-diagonal}
\end{figure}

The three closed-loop experiments with the corner-augmented surrogate
are repeated with the tightened upper bound of
Eq.~\eqref{eq:tightening}, i.e.\ the optimizer enforces
$\widehat{T}_{\max} \le 787$~K while all statistics continue to be
scored against the true process bound of $800$~K; controller, weights,
and surrogate are otherwise unchanged. The bottom row group of
Table~\ref{tab:closedloop} summarizes the outcome: the true peak
temperature never exceeds $800$~K on any of the three scan paths. The
sporadic single-step exceedances of the un-tightened corner-augmented
controller are eliminated on the seen vertical path ($805.8 \to
796.3$~K maximum) and spiral path ($802.7 \to 799.5$~K), as well as on
the unseen diagonal ($801.4 \to 791.8$~K,
Fig.~\ref{fig:cl-tight-diagonal}).
 
The spiral run illustrates the margin operating at its design limit:
the ground-truth maximum of $799.5$~K exceeds the \emph{tightened}
bound by $12.5$~K, an instantaneous under-prediction just inside the
calibrated $\delta = 13$~K, yet remains below the true $800$~K limit.
The worst-case surrogate error observed in these runs is thus
commensurate with, and covered by, the offline-calibrated quantile,
which is precisely the behavior the calibration is meant to deliver; a
larger error remains possible in principle, since the $95$th-percentile
margin is not a worst-case bound.
 
The conservatism cost of the margin is small. Because the hard bound is
active only at the peaks of the closed-loop temperature oscillation,
tightening it by $13$~K clips those peaks without shifting the mean
operating level: relative to the un-tightened corner-augmented runs,
the mean applied power changes by at most $0.06$~W ($< 0.5\%$) and the
occupancy of the soft lower bound is essentially unchanged (vertical:
$14.8\% \to 14.8\%$; spiral: $23.2\% \to 23.2\%$; diagonal: $35.8\% \to
36.8\%$ of settled steps). The main side effect appears in the
actuation activity on the diagonal path, where the mean step-to-step
power change grows from $0.58$ to $0.79$~W as the controller works
against the lower ceiling in the vicinity of the corner; the rate
remains well within actuator-realistic limits ($< 4\%$ of $P_{\max}$
per step). In combination, the corner-rich training data and the
calibrated margin thus deliver strict satisfaction of the true thermal
limit on both seen and geometrically unseen scan paths, at a
performance cost that is negligible for the seen geometries and modest
for the unseen one.
 
\subsection{Ablations and baseline comparisons}
\label{sec:ablations}
 
\begin{figure}[tb]
    \centering
    \includegraphics[width=0.95\linewidth]{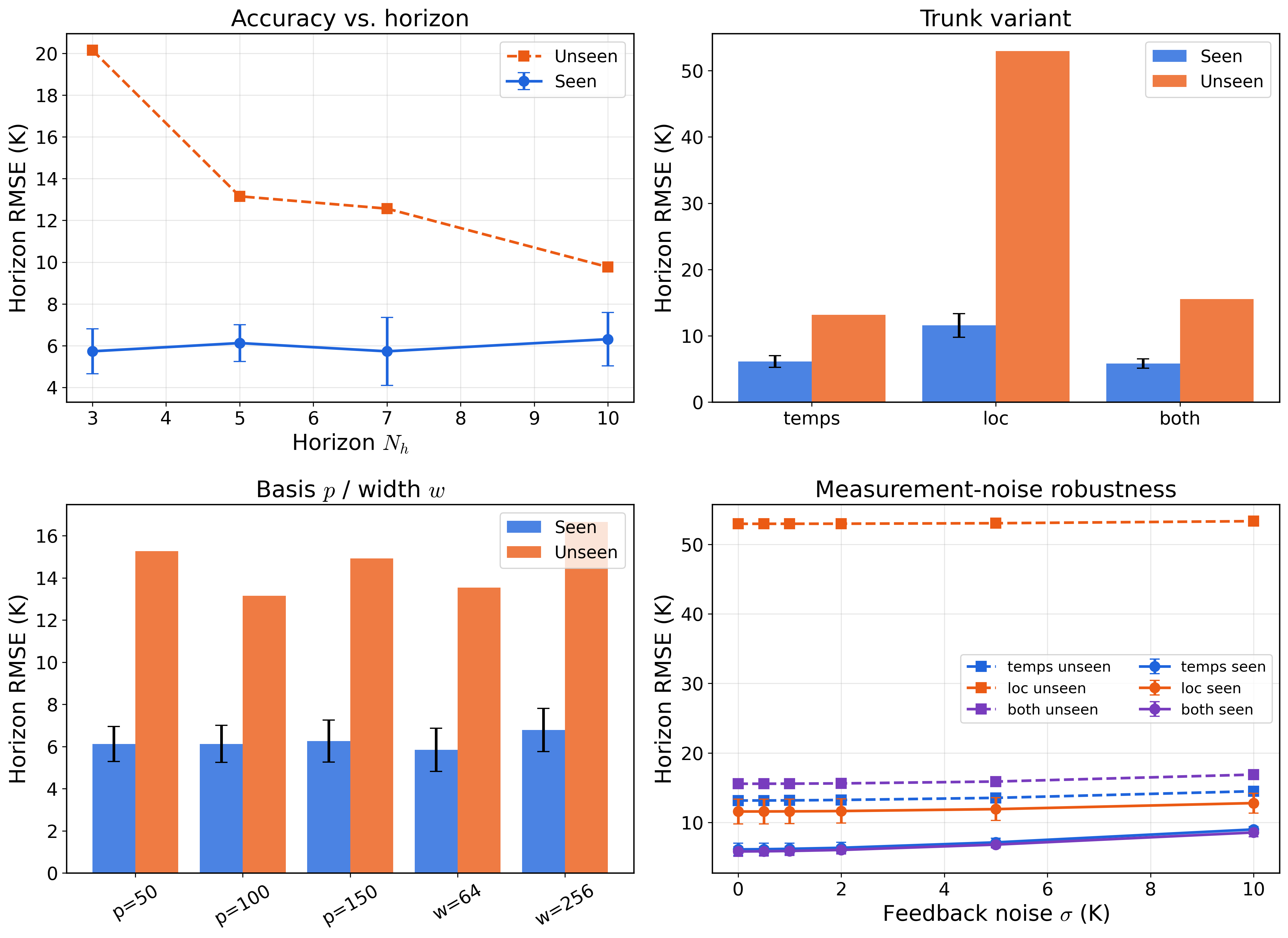}
    \caption{Offline ablation summary (blind-protocol horizon RMSE;
    seen = mean over horizontal/vertical/spiral with error bars, unseen = diagonal): accuracy vs.\ horizon $N_h$; trunk-input variant
    (temps / loc / both); basis dimension $p$ and hidden width $w$; and
    degradation under Gaussian noise on the temperature feedback.}
    \label{fig:ablation-offline}
\end{figure}

Guided by the finding of Section~\ref{sec:offline-results} that
aggregate open-loop accuracy does not reliably predict closed-loop
behavior, we evaluate each ablation axis at the level at which its
outcome actually matters. 
Architecture and robustness questions:
trunk-input design, basis dimension and network width, and
measurement-noise sensitivity, are assessed offline, where the
ranking is trustworthy and the mechanisms are interpretable
(Fig.~\ref{fig:ablation-offline}). 
The two axes that decide the
deployed controller, the horizon length $N_h$ and the one-shot
versus auto-regressive prediction structure, are evaluated in the
closed loop itself (Table~\ref{tab:ablation-cl}). 
All offline numbers
below are horizon RMSE values under the protocol of
Section~\ref{sec:offline-validation}, reported as the mean over the three
seen geometries (horizontal, vertical, spiral) and separately for the
unseen diagonal; the ablation models are independent retrains on
$\mathcal{D}_{\mathrm{corner}}$ with a fixed seed, so their absolute
values differ slightly from the deployed paper model.
 
\paragraph{Trunk-input variant (offline)}
We compare the paper's trunk encoding (\emph{temps}: $T_{\max}(t_k)$
plus the temperatures at the $N_h$ future laser locations) against a
variant that replaces the spatial pre-heat preview by the coordinates
of the current peak (\emph{loc}: $T_{\max}, X_{\max}, Y_{\max}$) and a
combination of both (\emph{both}). The result is clear: removing
the look-ahead temperatures collapses the surrogate, with the seen-path
RMSE nearly doubling ($6.1 \to 11.6$~K) and the unseen-diagonal RMSE
quadrupling ($13.2 \to 53.0$~K), whereas adding the peak location on
top of the temperatures yields no improvement (\emph{both}: $5.8$~K
seen, $15.6$~K unseen). The noise study confirms the mechanism: the
\emph{loc} variant's error is almost completely insensitive to noise
injected into its temperature inputs, confirming that without the
spatial preview the network has little usable thermal state to draw on.
The look-ahead temperatures of Eq.~\eqref{eq:trunk} are thus the
essential carrier of the pre-heating information that governs
corner accumulation, empirically validating the input design of
Section~\ref{sec:inputs}.
 
\paragraph{Basis dimension and network width (offline)}
Varying the number of basis functions $p \in \{50, 100, 150\}$ and the
hidden width $w \in \{64, 128, 256\}$ around the deployed configuration
($p=100$, $w=128$) changes the seen-path RMSE by less than $1$~K
($5.8$--$6.8$~K) and the unseen-diagonal RMSE by less than $3.5$~K,
with no systematic trend; the largest model ($w=256$) is in fact
slightly worse, consistent with mild overfitting at fixed data budget.
The accuracy is therefore insensitive to the network size over the
tested range, and the compact deployed configuration is retained for
its lower NLP cost.
 
\paragraph{Measurement-noise robustness (offline)}
Gaussian noise of standard deviation $\sigma$ is injected into all
temperature-derived trunk inputs, emulating thermal-camera measurement
error. The surrogate degrades slowly: for the deployed encoding the
seen-path RMSE grows from $6.1$~K (noise-free) to $6.4$~K at
$\sigma = 2$~K and $9.0$~K at $\sigma = 10$~K, while the
unseen-diagonal RMSE is nearly unaffected ($13.2 \to 14.5$~K). Even at
noise levels comparable to the surrogate's own error, the prediction
quality remains adequate for the constraint-tightened controller, and
the sensitivity is considerably less than one (an input noise of $10$~K induces
less than $3$~K of additional output error).
 
\paragraph{Horizon length (offline vs.\ closed loop)}
The horizon axis most clearly illustrates why closed-loop evaluation is essential. Offline, the blind-protocol RMSE on the unseen diagonal
improves monotonically with the horizon, from $20.2$~K at $N_h = 3$ to
$9.8$~K at $N_h = 10$ (Fig.~\ref{fig:ablation-offline}) -- partly a
protocol effect, since longer windows are re-anchored to the ground
truth less frequently per predicted step. In closed loop
(Table~\ref{tab:ablation-cl}, top group) this ordering disappears: the
one-step residual that the controller actually acts on saturates beyond
$N_h \approx 5$ (vertical: $2.9$--$3.6$~K across all horizons;
diagonal: $11.5$--$17.4$~K), the sporadic hard-bound exceedances are
small and non-monotone in $N_h$ (the largest, $7.1$~K, occurs at
$N_h = 7$), and the mean NLP solve time roughly triples from $29$~ms at
$N_h = 3$ to $90$~ms at $N_h = 10$ as the network evaluations and
decision variables grow with the horizon. The offline-optimal
$N_h = 10$ is thus the closed-loop-most-expensive setting with no
measurable control benefit, whereas $N_h = 5$ balances preview length,
one-step accuracy, and solve time, and is retained for the deployed
controller.
 
\begin{table}[tb]
\centering
\caption{Closed-loop ablation and baseline results
($N_{\mathrm{sim}} = 320$ steps, no constraint tightening; statistics
vs.\ FDM ground truth over the settled window $k \ge 10$).  RMSE$_1$ refers to the one-step prediction residual; solve is the mean
IPOPT time per control step.}
\label{tab:ablation-cl}
\scriptsize
\begin{tabular}{llcccccc}
\toprule
Setting & Scan path & $\max T_{\max}$ (K) & Over.\ (K) & $n_{>800}$ & $\%_{<760}$ & RMSE$_1$ (K) & Solve (ms) \\
\midrule
$N_h=3$  & Vertical  & 796.0 & 0.0 & 0 & 18.7 & 2.9  & 29 \\
$N_h=3$  & Diagonal  & 803.6 & 3.6 & 2 & 47.7 & 17.4 & 30 \\
$N_h=5$  & Vertical  & 802.7 & 2.7 & 1 & 10.6 & 3.6  & 40 \\
$N_h=5$  & Diagonal  & 798.0 & 0.0 & 0 & 27.1 & 12.5 & 50 \\
$N_h=7$  & Vertical  & 801.8 & 1.8 & 1 & 21.0 & 3.3  & 51 \\
$N_h=7$  & Diagonal  & 807.1 & 7.1 & 1 & 34.2 & 11.5 & 51 \\
$N_h=10$ & Vertical  & 794.4 & 0.0 & 0 & 22.6 & 3.2  & 90 \\
$N_h=10$ & Diagonal  & 791.6 & 0.0 & 0 & 33.9 & 12.4 & 91 \\
\midrule
SS-AR    & Vertical  & 791.0 & 0.0 & 0 & 33.5 & 12.1 & 50 \\
SS-AR    & Spiral    & 798.5 & 0.0 & 0 & 40.6 & 14.5 & 56 \\
SS-AR    & Diagonal  & 798.4 & 0.0 & 0 & 64.2 & 29.3 & 52 \\
\bottomrule
\end{tabular}
\end{table}
 
\begin{figure}[tb]
    \centering
    \includegraphics[width=0.95\linewidth]{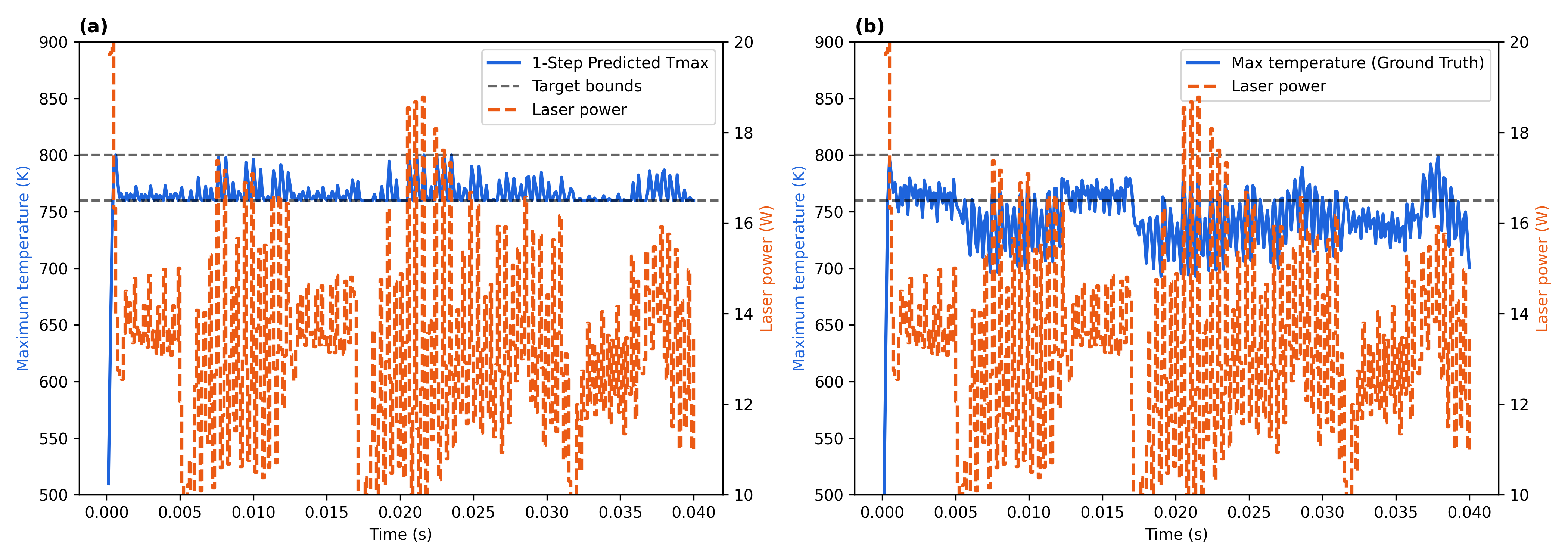}
    \caption{Closed-loop MPC on the unseen diagonal scan path with the
    single-step autoregressive (SS-AR) DeepONet baseline:
    (a)~surrogate one-step prediction and optimal power; (b)~FDM ground
    truth. The over-predicting recursive surrogate avoids the upper
    bound only by systematically under-powering: the true temperature
    spends $64\%$ of the settled steps below the $760$~K lower bound,
    dropping to $693$~K.}
    \label{fig:cl-ssar-diagonal}
\end{figure}

\paragraph{One-shot vs.\ autoregressive prediction (closed loop)}
As the structural baseline we deploy a single-step DeepONet in
autoregressive (SS-AR) mode: a network of identical architecture,
trained on the same corner-augmented data with the same branch inputs
$[P, X_c, Y_c, \dot{X}_c, \dot{Y}_c]$, predicts one step ahead and is
unrolled $N_h = 5$ times symbolically inside the NLP, each prediction
feeding the next step's trunk. Two structural consequences follow. The
first is recursive error compounding: offline, the SS-AR blind-window
error grows monotonically along the horizon (from $29$~K at the first
to $80$~K at the fifth step, vs.\ a flat $\approx 15$~K for the
one-shot model on the same data). The second is informational: under
autoregression the future thermal field is unavailable, so the SS-AR
trunk necessarily reduces to the scalar $T_{\max}(t_k)$ and the spatial
pre-heat preview, shown above to be the essential input, cannot be
supplied. The closed-loop consequences are shown in the bottom group of
Table~\ref{tab:ablation-cl} and Fig.~\ref{fig:cl-ssar-diagonal}. The
SS-AR controller avoids upper-bound violations, but it does so through
uncontrolled conservatism rather than accuracy: its one-step residual
is $3$--$4\times$ that of the one-shot model (e.g.\ $12.1$ vs.\ $3.6$~K
on the vertical path, $29.3$ vs.\ $14.8$~K on the diagonal) with a
strong over-prediction bias ($\bar{e} = +19.4$~K on the diagonal), so
the controller systematically under-powers: the plate spends $64\%$ of
the settled steps below the soft lower bound on the diagonal path,
dropping as far as $693$~K, and the actuation is $2$--$3\times$ rougher
($\overline{|\Delta P|} = 1.9$ vs.\ $0.6$~W). The process objective,
holding the temperature \emph{within} the window, is thus not met.
For completeness, a one-shot LSTM trained on identical data matches the
one-shot MS-DeepONet in offline accuracy but incurs roughly $6\times$
the inference latency and, owing to its recurrent gating, does not
admit the smooth algebraic embedding of
Section~\ref{sec:embedding}; it was therefore not deployed in
closed loop. The one-shot operator formulation thus combines
horizon-stable accuracy, the richest admissible state input, and the
most tractable NLP integration among the evaluated alternatives.
 
\subsection{Computational performance}
\label{sec:performance}

All timings were measured on a standard 14-core Intel machine (CPU
only; PyTorch~2.7, CasADi~3.7 with IPOPT), i.e.\ they reflect commodity hardware without any
GPU acceleration, code generation, or solver specialization.
We first compare the surrogate against the physics model over the
identical prediction task: one $N_h = 5$-step forecast of $\Tmax$ from
a representative operating state. A single MS-DeepONet forward pass
takes $0.42$~ms, whereas the five corresponding explicit FDM time steps
take $515$~ms, which is a speedup of roughly $1.2 \times 10^{3}$, comparable
to the three-orders-of-magnitude acceleration reported for the regional
Koopman surrogate on the same problem class \cite{zhang2026regional}. Embedded in the
receding-horizon optimization, one full MPC step -- an IPOPT solve
requiring on average $16$ iterations (at most $22$), each of which
evaluates and differentiates the surrogate graph -- takes $41$~ms on
average ($59.6$~ms at the 95th percentile, $71.7$~ms worst case),
consistent with the mean per-step solve times of $46$--$49$~ms logged
across the full closed-loop experiments of
Section~\ref{sec:results-cl}. A complete 320-step closed-loop run,
including the in-the-loop FDM plant, completes in about $50$~s.

These numbers place the approach precisely on the real-time map of this
process. The physical control interval is $\Delta t = 0.125$~ms, so the
FDM is $\sim\!4 \times 10^{3}$ times too slow even to \emph{simulate}
one horizon within one interval, whereas the surrogate forward pass
($3.4 \times \Delta t$) comes within a small factor of it.
The surrogate thus closes three of the roughly four orders of magnitude
separating the physics model from real time. The remaining gap is
introduced by the online optimization itself: at $\sim\!330 \times
\Delta t$, one IPOPT solve spans about $331$ control intervals, so the
present controller is a rigorous offline design-and-validation tool
rather than a real-time implementation. 

Two routes toward genuine real-time deployment at this time scale
follow directly from the structure of the method. First, the one-shot
multi-step prediction naturally supports a reduced re-solve rate: since
each solve returns a consistent $N_h$-step power plan, the controller
can apply several elements of the plan open-loop and re-solve only
every $m > 1$ intervals (move blocking), trading feedback bandwidth for
computation. 
Combined with code-generated solvers, warm-started SQP
iterations, and longer planning horizons, this shrinks the effective
per-interval cost substantially. Second, and more fundamentally, the
optimization itself can be amortized offline: a learning-to-optimize
surrogate trained to map the measured thermal state directly to the
optimal power sequence would reduce the online computation to a single
network evaluation, which, at $0.42$~ms, is already within a factor
of a few of the control interval on a plain CPU. We return to this
direction in Section~\ref{sec:discussion}.

\section{Discussion}
\label{sec:discussion}
\paragraph{One-shot versus autoregressive prediction}
The one-shot MS-DeepONet avoids the error compounding of a single-step autoregressive (SS-AR)
rollout. This is consistent with the bias/variance theory of multi-step forecasting: recursive
strategies accumulate error along the horizon, whereas direct or joint strategies avoid
accumulation at some variance cost~\cite{bentaieb2016bias}. Our SS-AR baseline, trained on
identical data and branch inputs, compounds from $29$ to $80$~K along the horizon and holds the
upper bound in closed loop only through a large $+19.4$~K over-prediction bias, pushing $64\%$ of
settled steps below the lower bound and roughening the actuation. The one-shot predictor incurs
no such trade-off. Autoregressive operation could nonetheless become preferable if the surrogate
outputs were made location-aware, predicting $(X_{\max},Y_{\max})$ alongside $T_{\max}$, so that
the rollout tracks the moving peak rather than a fixed observable; this is a promising route when
horizons longer than those studied here are required.

\paragraph{Relation to convex Koopman surrogate control}
The present operator-network controller and the regional Koopman controller of
\cite{zhang2026regional} occupy complementary points on an expressivity--tractability--generalization trade-off for the same problem class. The Koopman surrogate is linear in a lifted
state, yields a convex quadratic program, and generalizes by construction, at the cost of limited
expressivity. The operator network is more expressive and captures the strongly nonlinear,
history-dependent corner behavior, at the cost of a nonconvex nonlinear program and a
generalization that depends on training-data coverage. We make no quantitative head-to-head
benchmark; rather, the two are complementary components of a digital-twin workflow, to be selected
by the required balance of expressivity, online tractability, and guaranteed generalization.

\paragraph{Insufficiency of open-loop validation}
Two surrogates with nearly identical aggregate offline error on the unseen geometry produced
drastically different closed-loop outcomes. The reason is structural: the failure regime, the
source re-entering pre-heated material near the upper operating bound, is rarely visited by random
open-loop excitation, yet is permanently occupied by a controller that operates at the constraint.
Certifying a learned surrogate for control therefore requires closed-loop, regime-specific
evaluation; pooled open-loop metrics can mask a control-unsafe surrogate.

\paragraph{Limitations}
Generalization is data-coverage dependent, as the two-ensemble study makes explicit. The
constraint margin $\delta$ is an empirical quantile of the one-step under-prediction residual and
carries no formal probabilistic guarantee, in contrast to tube-based tightening for linear systems
with bounded disturbances~\cite{mayne2005robust}. Feedback is simulated by the finite-difference
plant emulating a thermal camera rather than measured. Actuation is power-only along fixed paths,
and the scope is single-track, single-layer.

\paragraph{Routes to real time}
The controller is presented as an offline design-and-validation tool: the mean solve time of
$41$~ms is roughly $330$ times the control interval, matching the operating mode of the companion
Koopman study~\cite{zhang2026regional}. The measured $0.42$~ms surrogate forward pass, about
$1.2\times10^{3}$ times faster than the equivalent finite-difference steps and $3.4$ times the
control interval, nonetheless enables two concrete paths to real-time operation. First, reducing
the re-solve rate through move blocking~\cite{cagienard2007moveblocking}: the one-shot $N_h$-step
plan can be applied open-loop for several steps between solves, combined with code-generated
warm-started solvers~\cite{andersson2019casadi,wachter2006ipopt}. Second, amortized
optimization (learning-to-optimize): a network trained to map the state directly to the optimal
power sequence would reduce the online cost to a single forward pass.

\section{Conclusion}
\label{sec:conclusion}
This paper pursued a twofold objective: to develop and validate a complete
surrogate-based pipeline for regulating the peak surface temperature of a
moving laser source, and to determine under what conditions the learned
surrogate can be trusted inside the control loop. 
Toward the first objective,
we presented a moving-source-aware multi-step deep operator network embedded
as a smooth nonlinear program inside a receding-horizon controller, holding
the peak temperature within a $[760, 800]$~K process window.
Against a high-fidelity finite-difference plant,
the surrogate forward pass is about $1.2\times10^{3}$ times faster than the equivalent solver steps.
Toward the
second objective, we exposed and diagnosed an out-of-distribution failure
mode, a 91~K corner under-prediction at a sharp path reversal on a
geometrically unseen path, and mitigated it: a controlled two-ensemble data
design reduced the failure to 1.4~K, and a calibrated one-sided constraint
margin of 13~K, set from the 95th percentile of the held-out one-step
under-prediction residual, achieved zero violations of the true 800~K bound
on all tested paths at negligible conservatism cost. 
The central
methodological finding, answering the second objective, is that aggregate
open-loop accuracy did not distinguish a control-ready surrogate from a
control-unsafe one; only closed-loop, regime-specific evaluation revealed
the difference.

Future work follows directly from the limitations. Location-aware outputs $(X_{\max},Y_{\max})$ would
enable stable autoregressive operation over longer horizons. Uncertainty quantification and
robust or quantile MPC would replace the empirical margin with a controller carrying explicit
guarantees. Joint power--velocity actuation and non-fixed paths would enlarge the admissible set.
Experimental validation with a real thermal camera, and extension to multi-track and multi-layer
builds, would close the gap to deployment. Finally, amortized optimization would convert the offline
design tool into a real-time controller.


\section{Conflict of interest}
The authors declare no conflicts of interest. 

\section*{Data availability}
The code and data supporting this study will be made available in a public repository upon publication of the final version of the paper.

\section*{Acknowledgments}
This research was supported by TUM International Graduate School of Science and Engineering (IGSSE) within the scope of the International Project Team OptiEnv. This research was funded in part by the Air Force Office of Scientific Research under grant FA9550-21-1-0381 and by the Office of Advanced Scientific Computing Research (ASCR) within the Department of Energy Office of Science under award number DE-SC0023163.

\bibliographystyle{elsarticle-num}
\bibliography{references}

\appendix
\section{Closed-loop results for the spiral scan path}
\label{app:spiral}
 
\begin{figure}[bth]
    \centering
    \includegraphics[width=0.95\linewidth]{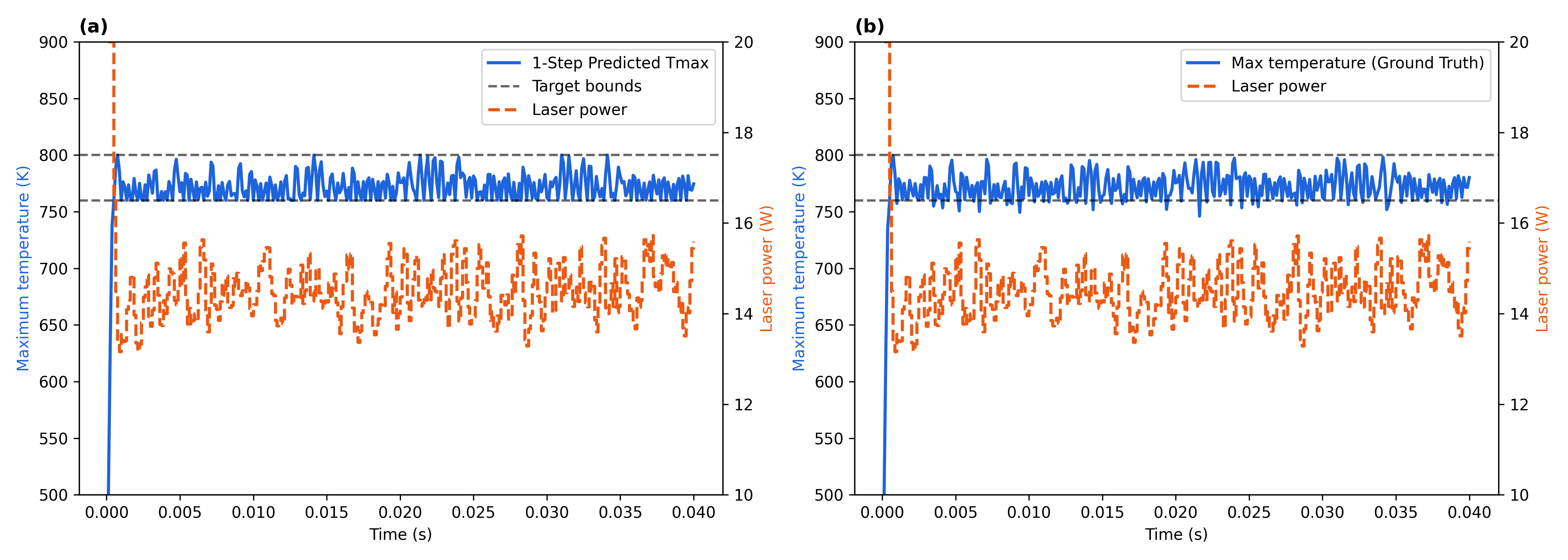}
    \caption{Closed-loop MPC on the seen spiral scan path with the
    baseline surrogate: (a)~surrogate prediction and optimal power;
    (b)~FDM ground truth. No hard-bound violation occurs.}
    \label{fig:app-spiral-base}
\end{figure}

\begin{figure}[bth]
    \centering
    \includegraphics[width=0.95\linewidth]{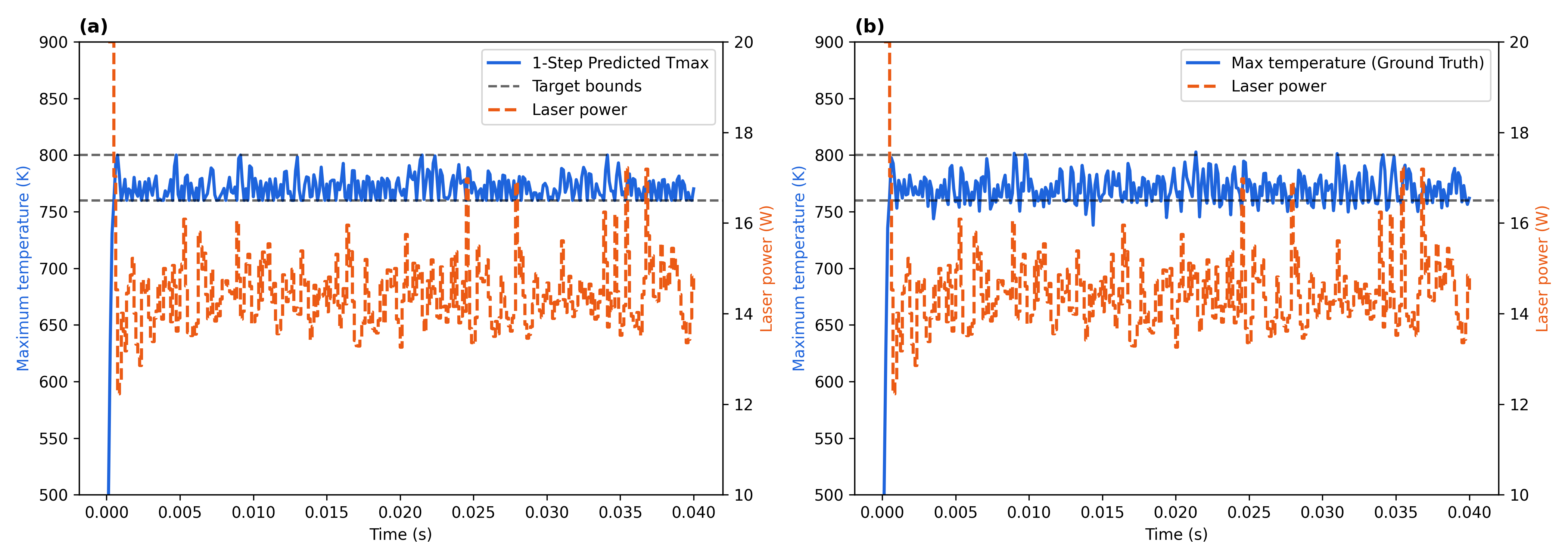}
    \caption{Closed-loop MPC on the seen spiral scan path with the
    corner-augmented surrogate. Isolated single-step exceedances of the
    hard bound of up to $2.7$~K reflect the increased in-distribution
    residual of the augmented model.}
    \label{fig:app-spiral-corneraug}
\end{figure}

This appendix collects the closed-loop results for the seen spiral scan
path, which corroborate the findings reported for the vertical path in
Sections~\ref{sec:results-cl-seen} and~\ref{sec:results-cl-ood}; the
corresponding summary statistics are included in
Table~\ref{tab:closedloop}. With the baseline surrogate
(Fig.~\ref{fig:app-spiral-base}) the controller regulates $\Tmax$
within the band without any hard-bound violation (maximum $798.3$~K);
the closed-loop residual (RMS $3.6$~K, maximum $14.0$~K) is moderately
larger than on the vertical path, consistent with the offline trend of
Section~\ref{sec:offline-results}, and essentially unbiased
($\bar{e} = 0.0$~K). With the corner-augmented surrogate
(Fig.~\ref{fig:app-spiral-corneraug}) the in-distribution residual
grows (RMS $6.9$~K, maximum $23.9$~K), and the ground truth exceeds the
hard bound in five isolated steps by at most $2.7$~K -- the same
trade-off observed on the vertical path.
 
\begin{figure}[t!]
    \centering
    \includegraphics[width=0.95\linewidth]{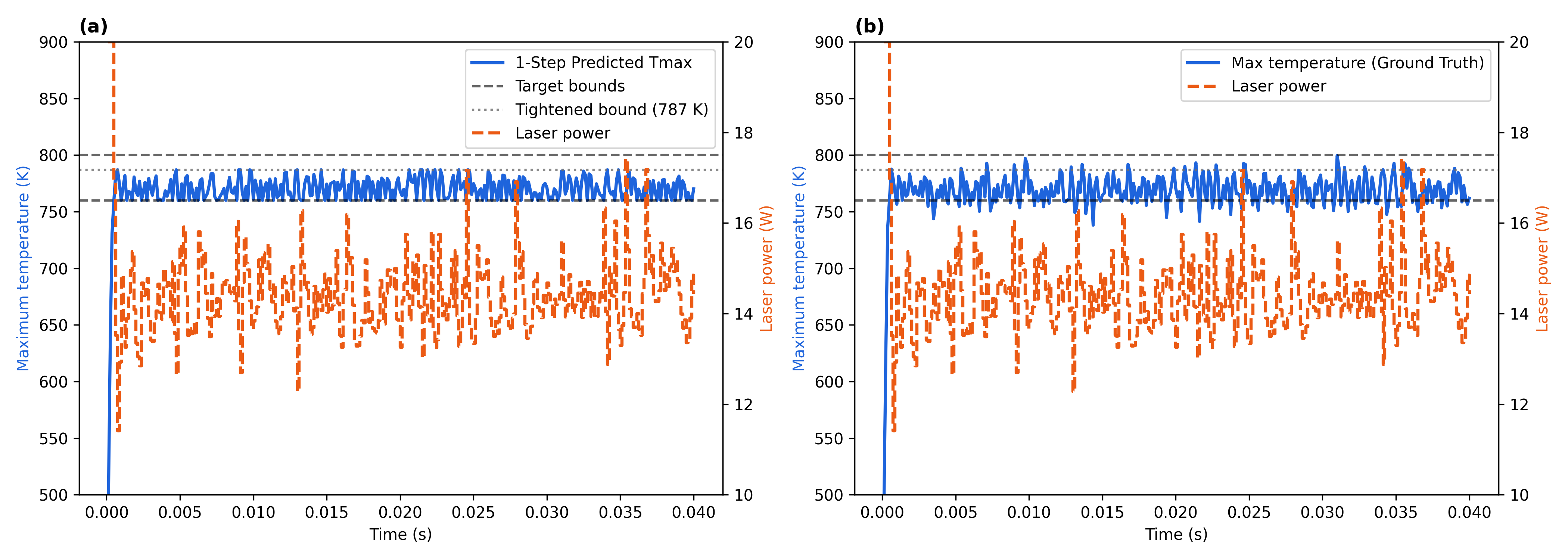}
    \caption{Closed-loop MPC on the seen spiral scan path with the
    corner-augmented surrogate and the tightened upper bound
    ($787$~K, dotted line). The true peak temperature remains below
    $800$~K throughout ($\max\, 799.5$~K, within the calibrated
    margin).}
    \label{fig:app-spiral-tight}
\end{figure}

With the tightened upper bound of Eq.~\eqref{eq:tightening} (Figure~\ref{fig:app-spiral-tight}) the five single-step exceedances of
the un-tightened run are eliminated; the ground-truth maximum of
$799.5$~K exceeds the tightened $787$~K bound by $12.5$~K -- just
inside the calibrated $\delta = 13$~K -- while remaining below the true
$800$~K limit, and the mean applied power is essentially unchanged
($14.40 \to 14.36$~W).

\end{document}